\documentclass[12pt,a4paper]{article}
\usepackage{amsmath,mathtools,pifont,amsfonts,amssymb,latexsym,multicol}
\usepackage{float,xspace,calc,theorem,ifthen}
\usepackage{fancyhdr,enumerate,epsfig}
\usepackage{times,epic,amscd} 
\usepackage{version,xr}
\usepackage{hyperref}
\usepackage{rotating}
\usepackage{color}
\usepackage{graphicx}
\newcommand{\mathpath}

\newcommand{\extension}
{pdf}

\hypersetup{colorlinks,pdfauthor={Andreas Knauf}}

\numberwithin{equation}{section}
\numberwithin{figure}{section}

\newcounter{antw}[section]

\theoremstyle{change}
\newtheorem{theorem}{Theorem} [section]
\newtheorem {lemma}[theorem]{Lemma}
\newtheorem {prop}[theorem]{Proposition}
\newtheorem {definition}[theorem]{Definition}

{\theorembodyfont{\normalfont}\newtheorem {remark}[theorem]{Remark}}
{\theorembodyfont{\normalfont}}
{\theorembodyfont{\normalfont}\newtheorem {example}[theorem]{Example}}
\newcommand{\beq}{\begin{equation}}
\newcommand{\eeq}{\end{equation}}
\newcommand{\Leq}[1]{\label{#1}\end{equation}}
\newcommand{\beqn}{\begin{eqnarray}}
\newcommand{\eeqn}{\end{eqnarray}}
\newcommand{\beqno}{\begin{eqnarray*}}
\newcommand{\eeqno}{\end{eqnarray*}}

\renewcommand {\l}{\left}
\newcommand {\ri}{\right}
\newcommand {\vep}{\varepsilon}

\newcommand {\LA}{\left\langle}
\newcommand {\RA}{\right\rangle}
\newcommand {\pa}{\partial}
\newcommand {\eh}{{\textstyle \frac{1}{2}}}

\newcommand {\ar}{\rightarrow}

\newcommand {\bN}{{\mathbb N}}
\newcommand {\bR}{{\mathbb R}}

\newcommand{\idty}{{\rm 1\mskip-4mu l}} 
\newcommand{\cA}{{\cal A}} %
\newcommand{\cC}{{\cal C}} %
\newcommand{\cD}{{\cal D}} %
\newcommand{\cF}{{\cal F}} %
\newcommand{\cH}{{\cal H}}

\newcommand{\cL}{{\cal L}}
\newcommand{\cM}{{\cal M}}

\newcommand{\cO}{{\cal O}} 
\newcommand{\cP}{{\cal P}}
\newcommand{\cT}{{\cal T}}

\newcommand{\ov}{\overline}
\newcommand{\bem}{\l(\! \begin{array}}
\newcommand{\eem}{\end{array}\!\ri)}
\newcommand{\bsm}{\left(\begin{smallmatrix}} 
\newcommand{\esm}{\end{smallmatrix}\right)}  
\newcommand{\NN}{\nonumber}

\newcommand{\qmbox}[1]{\quad\mbox{#1}\quad}
\renewcommand {\max}{{{\rm max}}}

\newcommand{\De}{\Delta^E} 
\newcommand{\Di}{\Delta^I} 

\begin{document}
\title {Asymptotic velocity for four celestial bodies}
\author{Andreas Knauf\thanks{
Department of Mathematics,
Friedrich-Alexander-University Erlangen-N\"urnberg,
Cauerstr.\ 11, D-91058 Erlangen, 
Germany, \texttt{knauf@math.fau.de}}}
\date{June 21, 2018}
\maketitle
\begin{abstract}
Asymptotic velocity is defined as the Ces\`aro  limit of velocity.
As such, its existence has been proven for bounded interaction potentials. 
This is known to be wrong in celestial mechanics with four or more bodies. 

Here we show for a class of pair potentials including the homogeneous ones of degree $-\alpha$
for $\alpha\in(0,2)$, that  asymptotic velocities
exist for up to four bodies, dimension three or larger, for any energy 
and {\em almost all\,} initial conditions on the energy surface. 
\end{abstract}
%
%
\section{Introduction}
%
In classical scattering theory one considers the motion of $n$ particles of masses $m_i>0$
and positions $q_i$ in $d$ spatial dimensions, generated by the Hamiltonian function
\beq 
\textstyle
{H}(p,q) := K(p) + V(q)\quad,\quad \mbox{with } 
V(q) := \sum_{1\le i<j\le n}V_{i,j}(q_i-q_j) ,
\Leq{Ham}
and kinetic energy $K(p) := \sum_{i=1}^n \frac{\|p_i\|^2}{2m_i}$.
The potential $V$
is assumed to be twice continuously differentiable and
{\em long-ranged}, that is, for some $\vep>0$, $I>0$ and in multi-index notation with $\beta\in \bN_0^d$
\beq
|\pa_\beta V_{i,j}(q)| \le  I \|q\|^{-\vep - |\beta|} \qquad \big(|\beta|\le 2,\, \|q\|\in[1,\infty)\, ,1\le i<j\le n \big).
\Leq{long:range}
For long-ranged $V_{i,j}\in C^2(\bR^d,\bR)$ and phase space $P:= T^*\bR^{dn}$, 
the Hamiltonian flow exists for all times (the escape times equal $T^\pm\equiv \pm\infty$).
For all initial conditions the existence of the {\em asymptotic velocities}, that is, the Ces\`aro limits  
\beq
\ov{v}^\pm(x_0):=\lim_{t\to T^\pm(x_0)} \frac{q(t,x_0)}{t}\in \bR^{dn}\qquad (x_0\in P)
\Leq{asymptotic:velocity}
of velocity $v$ (with momenta $p_i=m_iv_i$), was proven in \cite[Theorem 3.1]{De}. 
So in the infinite past and future the particles {\em cluster}, {\em i.e.} form set partitions of $N:=\{1,\ldots,n\}$
with equivalence relations given by $i \sim j$ if $\ov{v}^\pm_i(x_0)=\ov{v}^\pm_j(x_0)$.

Existence of asymptotic velocities is a fundamental result 
in classical scattering theory, 
whose proof is quite involved (see also \cite[Chapter 12.6]{Kn}).
As the $\ov{v}^\pm$ are constants of the motion, they present a weak form of integrability.
In fact, in the non-empty open set of $x_0$ leading to a trivial cluster partition 
(all $\ov{v}^+(x_0)$ being different), $\ov{v}^+$ is smooth if \eqref{long:range} is valid for all multiindices
$\beta\in \bN_0^{nd}$. However, in the presence of nontrivial clusters, one would need more
constants of the motion or other asymptotic information \cite{MS1,Sa} to have 
the 
so-called {\em asymptotic completeness}, see \cite{De}.

For the unbounded potentials considered here,
the question of classical asymptotic completeness is even more involved.
In the case of celestial mechanics, that is 
(in units where the gravitational constant equals one), 
$V_{i,j}(q)= -\frac{m_im_j}{\|q\|}$,
there are collisions for $n\ge 2$, and non-collision singularities for $n\ge4$, see \cite{Xi,Ge,Xu} and
Section \ref{sect2}.
Moreover, as 
shown in \cite{SX} for $d=1$ and $n=4$
there are initial conditions which (after regularization of binary collisions) do not lead to singularities but
for which $\lim_{t\to+\infty} \|q(t,x_0)\|/t=\infty$, so that asymptotic velocity \eqref{asymptotic:velocity}
does not exist (see also \cite{SD}).

Conversely, for $n<4$ celestial bodies asymptotic velocities exist for all initial conditions.
This follows from von Zeipel's Theorem and the fact that $j^+(x)<\infty$ for $n\le3$, see
Theorem \ref{thm:von:Zeipel}.
\begin{definition}\label{admissible}
The potential $V$ in \eqref{Ham} is {\bf admissible} if it is long-ranged,
central ($\,V_{i,j}(q)=W_{i,j}(\|q\|)$ for $q\in \bR^d\!\setminus\!\{0\}$),
and for some $\alpha\in (0,2)$ and $I>0$
\beq
|\pa_\beta V_{i,j}(q)| \le  I \|q\|^{-\alpha - |\beta|} \qquad(|\beta|\le 2,\, \|q\|\in(0,1]\, ,1\le i<j\le n).
\Leq{def:moderate}
\end{definition}
\begin{remark}[The constant $I$]\label{rem:constant:I}
By taking the maximum, we will use the same constant $I$ in \eqref{long:range} and \eqref{def:moderate}.
To gain more flexibility, by further increasing $I$, we assume that  \eqref{def:moderate} is valid
even for $\|q\|\le 2(\sum_{i=1}^n m_i)/\min_i m_i$.
\hfill $\Diamond$
\end{remark}
\begin{example}[Homogeneous potentials are admissible] \quad\\ 
If the pair potentials are homogeneous, that is, for $\alpha_{i,j}\in (0,2)$ and $I_{i,j}\in\bR$
\[V_{i,j}(q) = I_{i,j}\|q\|^{-\alpha_{i,j}} \quad(q\in \bR^d\!\setminus\!\{0\}\, ,1\le i<j\le n),\]
then $V$ is admissible. This includes gravitational and electrostatic interactions.\,$\Diamond$
\end{example}
\medskip

\noindent
Our main result is the following.\\[1mm]
{\bf Theorem }\label{thm:main} 
{\em For $n=4$ bodies, $d\ge3$ dimensions, and  an admissible potential $V$
the sets $\operatorname{NAV} := \{x\in P\mid \ov{v}^+(x) \mbox{ does not exist}\, \}$ 
of {\bf no asymptotic velocity}\qmbox{and}
\[
\operatorname{NAV}_{\!E} := \operatorname{NAV}\cap\, \Sigma_E \qmbox{with}
\Sigma_E:=H^{-1}(E)\quad(E\in\bR)\]
are Borel sets, with $\operatorname{NAV}_{\!E}\subseteq \Sigma_E$ being of 
Liouville measure $\sigma_E(\operatorname{NAV}_{\!E})=0$.}
\medskip

So in this case the asymptotic velocities \eqref{asymptotic:velocity} exist almost everywhere.\newline
Theorem \ref{thm:NAV:wandering} of Section \ref{sect2} analyzes the set $\operatorname{NAV}$
for $n\in \bN$ particles in arbitrary dimension $d$, interacting via long-ranged potentials.
\medskip

\noindent
{\bf Content }\\
In Section \ref{sect2} we will introduce some notation and present the main ideas of the proof.
A nondeterministic kinematical model for the asymptotics of those orbits, whose asymptotic
velocity does not exist, is exhibited in Section~\ref{sect3}. That the dynamics is in fact well described
by that model, will be shown in Section~\ref{sect4}. 
The proof of our main result is based on a Poincar\'e section method devised in \cite{FK1}.
The corresponding estimates form the content of the final Section \ref{sect5}.
\medskip

\noindent
{\bf Acknowledgements }
I am grateful to Stefan Fleischer for many conversations about the subject of non-collision singularities.
The critical remarks of the anonymous referee were very helpful for improving the text.

\section{Almost sure existence of asymptotic velocity} \label{sect2}
%
\subsection{Some general notation}
In general, for a $C^1$--vector field $X:P\to TP$ on a manifold $P$, 
we denote its maximal flow by $\Phi$. Then
\beq
\Phi\in C^1(D,P) \mbox{ with domain } D:=\{(t,x)\in \bR\times P\mid t\in(T^-(x),T^+(x))\}
\Leq{def:D}  
and its {\em escape times\,} $T\equiv T^+:P\to (0,+\infty]$ and $T^-:P\to [-\infty,0)$. 
$T$ is lower semicontinuous in the obvious topology of $(0,\infty]:=(0,\infty)\cup\{\infty\}$. 

We mark the values of a phase space variable $V\in C^1(P,\bR^d)$ along the flow line 
$t\mapsto \Phi_t(x)$ by a tilde:
\[\widetilde{V}(t):=V\circ \Phi_t(x).\]
For one-sided limits we use the notations 
\[  \widetilde{V}(t^+):=\lim_{s\searrow t} \widetilde{V}(s) \qmbox{and} 
\widetilde{V}(t^-):=\lim_{s\nearrow t}\widetilde{V}(s). \]

\smallskip

\noindent
Here we consider $n$ particles of masses $m_i>0$ in the configuration spaces
$M_i:=\bR^d$. 
On the joint configuration space $M :=\bigoplus _{i=1}^nM_i$
of all particles we use the inner product
$\LA\cdot,\cdot\RA_\cM$
generated by the mass matrix\label{mass matrix}
$\cM:={\rm diag}(m_1,\ldots,m_n)\otimes \idty_d$,
\beq
\LA\cdot,\cdot\RA_A:M\times M\to\bR \qmbox{,} 
\LA q,q'\RA_A :=  \LA q,A q'\RA,
\Leq{inner:product}
denoting the bilinear form for the 
matrix $A\in {\rm Mat}(M,\bR)$, with
the canonical inner product $\LA \cdot,\cdot\RA$.
We set $\|q\|_A:=\LA q,q\RA_A^{1/2}$.

The {\em collision set} in the configuration space is given by
\beq
\Delta := \{q\in M\mid q_i=q_j\mbox{ for some }i\neq j\in N\}.
\Leq{coll:set}
We consider potentials $V\in C^2(\widehat{M},\bR)$ 
on the {\em non-collision configuration space}
$\widehat{M} := M\backslash\Delta$. 
The Hamiltonian function $H\in C^2(P,\bR)$ on the phase space 
$P:=T^*\widehat{M}$ is given by \eqref{Ham},
with pair potentials $V_{i,j}$. 

With the Euclidean gradient $\nabla$ on $\bR^d$
the Hamiltonian equations of \eqref{Ham} are
\[\textstyle
\dot{p}_i= \sum_{j\in N\setminus \{i\}} \nabla V_{i,j}(q_j-q_i)
\qmbox{,}
\dot{q}_i=\frac{p_i}{m_i}
\qquad (i\in N).\]
Using the natural symplectic form $\omega_0$ on the cotangent bundle
$P $, we write these as
$\dot{x}=X_H(x)$ for the Hamiltonian vector field
$X_H$ defined by ${\bf i}_{X_H}\omega_0=dH$.

Since $H\circ {\cal R} = H$ for the involution ${\cal R}:P\to P$, $(p,q)\mapsto (-p,q)$, 
the maximal flow $\Phi\in C^1\big(D, P \big)$ is {\em reversible}: 
${\cal R}\circ\Phi_t\circ {\cal R}=\Phi_{-t}$.
In particular $T^-(p,q)=-T^+(-p,q)$.
$\Phi$ restricts to the {\em energy surfaces} $\Sigma_E:=H^{-1}(E)$.
We write $\big(p(t,x),q(t,x)\big):= \Phi(t,x)$ 
for the momenta resp.\ positions at time $t$ starting at $x\in P $, 
and given an initial condition $x$, we even write 
$\big(\tilde{p}(t),\tilde{q}(t)\big)$.

\begin{definition}\label{def:sing:coll} 
The set of phase space points {\bf experiencing a singularity} is\\[-3mm]
\[ \operatorname{Sing} := \big\{ x \in P \mid T(x) < \infty \big\}\, , \]
whereas its subsets experiencing a {\bf (non-)\,collision singularity} are
\beq
\operatorname{Coll}:= \big\{ x\in \operatorname{Sing}\mid
\lim_{t\nearrow T(x)}q(t,x)\ \mbox{exists} \big\}\qmbox{and}
\operatorname{NC} := \operatorname{Sing} \setminus \operatorname{Coll}.
\Leq{cs}
The corresponding orbits in $\operatorname{Coll}$ resp.\ $\operatorname{NC}$
are called {\bf (non-)\,collision orbits}. 
\end{definition}
$\operatorname{Sing}\cap \operatorname{NAV} = \operatorname{NC}$, and 
$\operatorname{Sing}$, $\operatorname{Coll}$ and $\operatorname{NC}$ are Borel subsets of $P$.

For arbitrary $n\in \bN$, $d\ge2$ and a large class of pair interactions including the homogeneous ones
(see Saari \cite{Sa1}) the set of initial conditions leading to collisions are of Liouville measure zero 
for all total energies $E$, see \cite{FK2}. 
Strongly based on finiteness of the escape time, the same result concerning non-collision singularities
was proven for $n=4$ by Saari in \cite{Sa2} and by Fleischer in \cite{Fl}. 

To show almost sure existence of asymptotic velocities, we have to 
take into account the case of solutions of the initial value problem that exist for all times, but where 
\eqref{asymptotic:velocity} does not exist.
So we have to argue differently from \cite{Sa2,Fl}.

\smallskip

\subsection{Cluster decompositions}
%
This section can be skipped, as we
recall the notations introduced here in later sections.
We introduce standard notions for the set partitions of $N=\{1,\ldots,n\}$: 
\vspace*{-5mm}
%
\begin{definition}\quad\\[-6mm]
\label{def:cluster}
\begin{enumerate}[$\bullet$]
\item
A \textbf{set partition} or \textbf{cluster decomposition} of $N$ is a set $\cC:=\{C_1,\ldots,C_k\}$ of
\textbf{blocks} or \textbf{clusters} $\emptyset\neq C_\ell\subseteq N$ such that\vspace*{-2mm}
\[\textstyle
\bigcup_{\ell=1}^kC_\ell=N \qmbox{and}
C_\ell\cap C_m=\emptyset\ \mbox{ for } \ell\ne m\, .
\]
We denote by $\sim_\cC$ (or $\sim$, if there is no ambiguity) the equivalence relation on $ N$ induced by 
$\cC$; the corresponding equivalence classes are denoted by $[\cdot ]_\cC$.
\item
The \textbf{lattice of partitions} $\cP(N)$ is the set of cluster
decompositions $\cC$ of $N$, partially ordered by
\textbf{refinement}, i.e., 
\[
\cC=\{C_1,\ldots,C_k\}\preccurlyeq \{D_1,\ldots,D_\ell\}=\cD\, ,
\] 
if $C_m\subseteq D_{\pi(m)}$ for an appropriate surjection
$\pi:\{1,\ldots,k\}\to\{1,\ldots,\ell\}$. 

In this case, $\cC$ is called  \textbf{finer} than $\cD$ and $\cD$
\textbf{coarser} than $\cC$.\\
The unique finest and coarsest elements of $\cP(N)$ are
\[\cC_{\wedge} := \big\{ \{1\}, \ldots, \{n\} \big\} \qmbox{and} \cC_{\vee} := \{N\}=\big\{\{1,\ldots,n\}\big\},\]
respectively.
Two clusters $\cC,\cD\in\cP(N)$ are {\bf comparable} 
if $\cC\preceq\cD$ or $\cC\succeq\cD$.
\item
The  \textbf{rank} of $\cC\in\cP(N)$ is the number  $|\cC|$ of its blocks, 
and 
\beq \cP_k(N):=\{\cC\in \cP(N) \mid |\cC|=k\}\qquad(k\in N). \Leq{Pk}
\item
The  \textbf{join} of $\cC$ and $\cD\in\cP(N)$, denoted as  
$\cC\vee\cD$, is the finest cluster decomposition that is coarser than
both $\cC$ and $\cD$.
\end{enumerate}
\end{definition}

We use the following partitions to decompose configuration space $M$:\\
given a subset $\emptyset \neq C \subseteq  N$, we define the corresponding \emph{collision set} as
\begin{equation*}
\Delta_C^E := \left\{ q \in M \ | \ q_i = q_j \mbox{ if } i,j \in C \right\},
\end{equation*}
and for a cluster decomposition $\cC$ we define the \emph{$\cC$-collision subspace}
\begin{equation}
\label{eq:DefDeltaE}
\Delta_\cC^E := \left\{ q \in M \ | \ q_i = q_j \mbox{ if } [i]_\cC = [j]_\cC \right\} = \bigcap_{C\in\cC} \Delta_C^E \, .
\end{equation}
By $\Pi_C^E$ we denote the $\cM$-orthogonal projection onto the subspace $\Delta_C^E$, 
and we denote the complementary projection $\idty_\cM - \Pi_C^E$ by $\Pi_C^I$.
Accordingly, we denote the ({\em external}) projection onto 
$\Delta_\cC^E$ by $\Pi_\cC^E := \prod_{C\in \cC} \Pi_C^E$, 
and the ({\em internal}) projection by $\Pi_\cC^I = \idty_\cM - \Pi_\cC^E = \sum_{C\in \cC} \Pi_C^I$.
The image of $\Pi_C^I$ is given by
\[
\textstyle 
\Delta_C^I := 
\left\{ q \in M \ \left| \ \sum_{i\in C} m_i q_i = 0,\ \forall \, i \in N\!\setminus\! C: \, q_i = 0 \,   \right.\right\},
\]
the image of $\Pi_\cC^I$ is $\
\Delta_\cC^I := \left\{ q \in M \ \left| \ \forall\,C\in \cC: \sum_{i \in C} m_i q_i = 0  \right.\right\} = 
\bigoplus_{C\in \cC} \Delta_C^I$. 

In particular, $\Delta_{\cC_{\wedge}}^E=M$, see Def.\ \ref{def:cluster}. 
We get a $\cM$-orthogonal decomposition
\begin{equation}
\label{M:E:I}
\textstyle
M = \Delta_\cC^E \oplus \bigoplus_{C\in \cC} \Delta_C^I.
\end{equation}
For a nonempty subset $C\subseteq N$ we define the 
{\em cluster mass\,}, {\em cluster barycenter}  and {\em cluster momentum} of $C$ by
\beq
\textstyle
m_C := \sum_{j\in C} m_j \qmbox{,}
  q_C := \frac{1}{m_C} \sum_{j\in C} m_j q_j \qmbox{and} p_C:=\sum_{i\in C}p_i.
\Leq{mqp:C}  
In particular $m_N$ equals the 
{\em total mass} of the particle system.
Then for the partitions $\cC\in \cP(N)$ the $i$--th component of the 
cluster projection $q^E_\cC := \Pi^E_\cC(q)$ is given by the barycenter, respectively 
its distance from the barycenter
\begin{equation}
\label{cl:bar}
\big(q^E_\cC\big)_i = q_{[i]_\cC} \qmbox{,}
\big(q^I_\cC\big)_i = q_i-q_{[i]_\cC} \mbox{ for } q^I_\cC := \Pi^I_\cC(q) \qquad (i\in N).
\end{equation}
Join of partitions corresponds to intersection of collision subspaces:
\[\Delta_\cC^{E}\cap \Delta_\cD^{E} = \Delta_{\cC\vee\cD}^{E}
\qquad\bigl(\cC,\cD\in \cP(N)\bigr). \]
So for $\cC\in \cP(N)$, the mutually disjoint sets
\beq
\Xi_\cC^{(0)} \,:=\, \Delta_\cC^{E}\,\big\backslash \,
\mbox{{\normalsize $\bigcup\limits_{\cD\succneqq\cC}$}}
\Delta_\cD^{E} 
\,=\, \{ q \in M \mid q_i = q_j \mbox{ if and only if  }i \sim_\cC j \} \  
\Leq{Xi:C:0}
form a set partition of $M$, with $\Xi_{\cC_{\wedge}}^{(0)}=\widehat{M}$. 


Denoting by $M^*$ the dual space of the vector space $M$, there are natural 
identifications $TM\cong M\times M,\ T^*M\cong M^*\times M$ of the tangent
space resp.\ phase space of $M$. This gives rise to the inner products
\[\LA \cdot,\cdot\RA_{TM}:TM\times TM\to\bR \qmbox{,} 
\LA (q,v),(q',v')\RA_{TM}:=\LA q,q'\RA_\cM+\LA v,v'\RA_\cM\]
and
\[\LA \cdot,\cdot\RA_{T^*M}:T^*M\times T^*M\to\bR \qmbox{,} 
\LA (p,q),(q',p')\RA_{T^*M}:= \LA q,q'\RA_\cM+\LA p,p'\RA_{\cM^{-1}}\]
(with $\LA p,p'\RA_{\cM^{-1}} = \sum_{i=1}^n \LA p_i, p_i' \RA / m_i$  
for the momentum vector $p=(p_1,\ldots,p_n) )$.

The tangent space $TU$ of any subspace $U\subseteq M$ is naturally a 
subspace of $TM$. 
Using the inner product, we also consider $T^*U$ as a subspace of $T^*M$.

We thus obtain $T^*M$--orthogonal decompositions
\[\textstyle 
T^*M = T^*\De_\cC\oplus\bigoplus_{C\in\cC}T^*(\Di_C)\qquad \big(\cC\in\cP(N)\big)\]
of phase space. 
With  
$\widehat{\Pi}^I_\cC := \idty_{T^*M}-\widehat{\Pi}^E_\cC
=\sum_{C\in\cC}\widehat{\Pi}^I_C$
the $T^*M$--orthogonal projections 
$\widehat{\Pi}^E_\cC, \, \widehat{\Pi}^I_\cC: T^*M \to T^*M$
onto these subspaces are given by the {\em cluster coordinates}
\beq
\textstyle
(p^E,q^E) := \widehat{\Pi}^E_\cC(p,q) \qmbox{with} (p_i,q_i)
= \big(\frac{m_i}{m_{[i]}}p_{[i]} ,  q_{[i]}\big) \quad(i\in N),
\Leq{cl:co}
and {\em relative coordinates}
\[(p^I,q^I) := \widehat{\Pi}^I_\cC(p,q) \qmbox{with} 
(p^I_i, q^I_i) = (p_i-p^E_i,q_i-q^E_i) \quad (i\in N).\]
Compared to (\ref{cl:bar}) we omitted the subindex $\cC$ in (\ref{cl:co}).

For nontrivial partitions 
neither the cluster coordinates nor the relative coordinates are
coordinates in the strict sense.
Later, however, we need such coordinates
on the above-mentioned symplectic subspaces of phase space.

\begin{lemma}\label{lem:sym}
The vector space automorphisms 
\beq
\big(\widehat{\Pi}^E_\cC,\widehat{\Pi}^I_\cC \big):T^*M\longrightarrow
T^*\De_\cC\textstyle \oplus\bigoplus_{C\in\cC}T^*(\Di_C)
\qquad\big(\cC\in \cP(N)\big)
\Leq{sy:trans}
are symplectic w.r.t.\ the natural symplectic forms on these cotangent bundles.
\end{lemma}
{\bf Proof.}
This follows from 
$T^*\big(\De_\cC\oplus\bigoplus_{C\in\cC}\Di_C\big)=T^*\De_\cC\oplus\bigoplus_{C\in\cC}T^*(\Di_C)$.\hfill$\Box$\\[2mm]
{\em Total angular momentum}
\beq
\textstyle
L:T^*M\to\bR^d\wedge\bR^d \qmbox{,} L(p,q) = \sum_{i=1}^nq_i\wedge p_i
\Leq{total:ang:momentum}
and {\em total kinetic energy} 
\[\textstyle
K:T^*M\to\bR \qmbox{,} K(p,q) \equiv K(p)=\eh\LA p,p\RA_{\cM^*} =
\sum_{i=1}^n\frac{\LA p_i,p_i\RA}{2m_i}\]
(using a sloppy notation) both split for $\cC\in \cP(N)$ into sums of {\em barycentric}
\beqn
L^E_\cC &:=& \textstyle
L\circ\widehat{\Pi}^E_\cC \qmbox{,} 
L^E_\cC(p,q) = \sum_{C\in\cC}q_C\wedge p_C,
\label{ex:ang:mom}\\
K^E_\cC &:=& \textstyle
K\circ\widehat{\Pi}^E_\cC \qmbox{,} K^E_\cC(p,q) = \sum_{C\in\cC}
\frac{\LA p_C,p_C\RA}{2m_C}\NN
\eeqn
and {\em relative} terms for the clusters $C\in\cC$
\beqn
L_C^I &:=& \textstyle
 L\circ\widehat{\Pi}_C^I \qmbox{,} 
 L_C^I(p,q) = \sum_{i\in C}q^I_i\wedge p^I_i
 \label{rel:ang:mom}\\
K_C^I &:=& \textstyle
K\circ\widehat{\Pi}_C^I \qmbox{,} 
 K_C^I(p,q) = \sum_{i\in C}\frac{\LA p^I_i,p^I_i\RA}{2m_i}.\NN
\eeqn
That is,
$L = L^E_\cC + \sum_{C\in\cC}L_C^I \qmbox{and} 
  K = K^E_\cC + K^I_\cC \qmbox{with}K^I_\cC:= \sum_{C\in\cC}K_C^I$.

The decomposition of the potential $V$ as a sum of external and internal terms is not given
by composition with the corresponding projections.\\
Instead, we set $V_C:= \sum_{i<j\in C}V_{i,j}$ for $C\subseteq N$ and
\beq
\textstyle
V^I_\cC:=\sum_{C\in \cC}V_C\qmbox{,}V^E_\cC:=V-V^I_\cC\qquad (\cC\in \cP(N)).
\Leq{V:IE:C}
Finally, the Hamiltonian function $H$ decomposes into the terms 
\beq
\textstyle
H^I_\cC := \sum_{C\in \cC}H^I_C \qmbox{with}
H^I_C := K^I_C + V^I_C\qquad (\cC\in \cP(N)),
\Leq{H:IE:C}
and similarly for $H^E_C$ and $H^E_\cC$.

For simplifying notation, we sometimes omit the super-index  $E$ (but not $I$).
%
\subsection{The wandering set}
%
We base our article on the observation that $\operatorname{NAV}$ is wandering, see 
Theorem \ref{thm:NAV:wandering}. 
We shortly discuss the general method from \cite{FK1} for $C^1$, volume-preserving dynamical systems
$(P,\Omega, X)$, with $X:P\to TP$ having Lie derivative $\cL_X\Omega=0$.\vspace*{-4mm}
\begin{definition}
The {\bf wandering set} $\operatorname{Wand} \subseteq P$ of the flow $\Phi$ 
consists of those $x\in P$ which have a neighborhood   $U_x$ 
so that for a suitable time $t_x\in(0,T^+(x))$
\[ U_x \cap \Phi \big( \big((t_x,T^+(x)) \times U_x\big) \cap D \big) = \emptyset.\]
\end{definition}
\smallskip

\noindent 
As $\cL_X\Omega=0$, $\Phi$ preserves the volume form $\Omega$. 
Consider now for the differential form of degree $\dim(P)-1$
\[{\cal V} := \boldsymbol{i}_{X} \Omega\]
a sequence of hypersurfaces 
${\iota}_m: \cH_m \subseteq P$ ($m\in\bN$) which fulfill the Assumptions
\begin{enumerate}
\item[\bf (A1)] 
to be transversal to $X$, so that ${\cal V}_m := {\iota}_m^*({\cal V})$ are volume forms on $\cH_m$,
\item[\bf (A2)]  
to have finite volumes ($\int_{\cH_m} {\cal V}_m < \infty$), and
$\lim_{m\to\infty}\int_{\cH_m} {\cal V}_m =0$.
\end{enumerate}
\beq
\operatorname{Trans}:=
\{x\in P \mid \exists\, m_0\in\bN\;\forall \,m\ge m_0: \cO^+(x)\cap \overline{\cH}_m\neq\emptyset \}
\Leq{trans:in:sing}
is the set of {\em transition points}, whose forward orbits eventually hit all of these surfaces.
Then
\begin{theorem}[\protect{\cite{FK1}, Theorem A}]\label{thm:FK} \
$\Omega(\operatorname{Wand}\cap\operatorname{Trans}) = 0$. 
\end{theorem}
So our task is to show that this {\em method of Poincar\'e surfaces} 
is applicable, {\em i.e.} setting $(P,\Omega,X) := (\Sigma_E, \sigma_E,X_H|_{\Sigma_E})$
for energy $E\in \bR$, (omitting indices $E$)
\begin{enumerate}[\bf (B1)]
\item 
to show that $\operatorname{NAV}\subseteq\operatorname{Wand}$, and
\item 
to define a sequence of hypersurfaces $\cH_m \subseteq \Sigma_E$, fulfilling Assumptions
{(A1)} and {(A2)}, for which
$\operatorname{NAV}\subseteq\operatorname{Trans}$.
\label{B2}
\end{enumerate}
%
\subsection{Proof of assertion (B1)} 
%
We consider the time evolution for the moment of inertia
\[\textstyle
J: P\to \bR\qmbox{,} J(q):= \eh\sum_{i=1}^n m_i \|q_i\|^2 = \eh \langle q,q \rangle_\cM\] 
and use the following generalization by Fleischer  \cite[Theorem 2.4.4]{Fl} of von Zeipel's theorem
\cite{Ze, Sp}. As it stands, it is a result about individual orbits.
\begin{theorem}[von Zeipel]\quad\\  
\label{thm:von:Zeipel}%
For $n\in\bN$ particles in $d$ dimensions, a long ranged potential $\,V$ (see \eqref{long:range}) and
\beq
\tilde{j}_x : \big( 0,T^+(x) \big) \to [0,\infty) \qmbox{,} \tilde{j}_x(t) := J\big( q(t,x)/t \big) \qquad (x\in P)
\Leq{j:x}
the limit $j^+(x):=\lim_{t\nearrow T^+(x)} \tilde{j}_x(t)$ exists in $[0,\infty]$.\,
$j^+(x)$ is finite if and only if the asymptotic velocity $\overline{v}^+(x)\in \bR^{dn}$ 
from \eqref{asymptotic:velocity} exists.
\end{theorem}
Indeed, in general the phase space function $j^+:P\to [0,\infty]$ is discontinuous.
\begin{example}[Discontinuity of the phase space function $j^+$]\quad\\[-6mm] 
\begin{enumerate}[1.]
\item 
For the Hamiltonian 
$H(p,q):= \eh p^2+V(q)$ with potential $V\in C^\infty_c(\bR,\bR)$ of compact support and a unique
absolute maximum $V(0)=1=:E$, $j^+$ is {\em not upper semicontinuous} at $x_0:=(0,0)\in \Sigma_E$.
Then $j^+(x_0)=0$, whereas
for all initial conditions $x'_0=(p'_0,q'_0)\in \Sigma_E$ with $p'_0 q'_0>0$ one has $j^+(x'_0)=\eh$. 
\item 
For the $d\ge2$--dimensional Kepler Hamiltonian 
$H(p,q):=\eh \|p\|^2-1/\|q-c\|$ with $c\in \bR^d\!\setminus\!\{0\}$, $j^+$ is {\em not lower semicontinuous} 
at initial conditions $x_0=(0,q_0)$ with $\|q_0-c\|=1$, that is, $E:=H(x_0)=-1$ 
and $x_0\in \operatorname{Coll}$.
Then $j^+(x_0)=\eh\|c/T^+(x_0)\|^2$, but for all 
$x'_0\in \Sigma_E\! \setminus \!\operatorname{Coll}$ one has $j^+(x'_0)=0$.
\hfill $\Diamond$
\end{enumerate}
\end{example}
However, we have the following result, proven in Section \ref{sub:sect:NAV:wanders},
which is of independent interest and proves assertion (B1). It uses the time dependent
{\em cluster function} $\cA:\big(0,T(x)\big) \to\cP(N)$, defined in Section \ref{sub:sect:NAV:wanders}, 
see \eqref{00V}.
\begin{theorem}[Initial conditions without asymptotic velocities]\quad \label{thm:NAV:wandering}\\
For $n\in\bN$ particles in $d\in \bN$ dimensions and a long-ranged potential $V$, 
\begin{enumerate}[1.]
\item 
For initial conditions $x\in \operatorname{NAV}$, the cluster-external speed has the limit\linebreak
$\lim_{t\nearrow T(x)}\|\frac{d}{dt}\tilde{q}_{\cA(t)}^E(t) \|_\cM=\infty$, so that
$\lim_{t\nearrow T(x)}V(q(t,x))=-\infty$.
\item 
 $j^+:P\to [0,\infty]$ is continuous\,%
 \footnote{regarding, as usual, $[0,\infty]$ as the one-point compactification of $[0,\infty)$} 
 at the points $x\in \operatorname{NAV}$.
\item 
Thus $\,\operatorname{NAV}\subseteq \operatorname{Wand}$, 
and $\,\operatorname{NAV}\subseteq P$ is a Borel set.
\end{enumerate}
\end{theorem}
%
%
\subsection{Proof of assertion (B2)} \label{sub:sect:B2}
%
Assertion (B2) on page \pageref{B2} states that there is a sequence of hypersurfaces 
$\cH_m\subseteq \Sigma_E$ , fulfilling Assumptions {(A1)} and {(A2)}, for which
$\operatorname{NAV}_E\subseteq\operatorname{Trans}_E$. 
We will perform two changes of that scheme, using the constants of motion $L$ and $p_N$.
First consider total momentum and the center of mass, in the notation \eqref{mqp:C}: 
\[\textstyle
p_N = \sum_{i\in N} p_i\qmbox{and}
q_{N} = \sum_{i\in N} \frac{m_i}{m_N} q_i\qmbox{, with}m_N=\sum_{i\in N} m_i \, .\]
Subsets of  the energy surface in the {\em center of mass frame}
\[\Sigma_{E,0}:=\{(p,q)\in \Sigma_E\mid p_N=0,\, q_N=0\}\]
are denoted by using the same indices $E,0$.
\begin{lemma}[Center of mass]  
If $\,\operatorname{NAV}_{\!E,0}\subseteq \Sigma_{E,0}$ are of Liouville measure zero for all $E\in \bR$,
then $\operatorname{NAV}_{E}\subseteq \Sigma_{E}$ are of Liouville measure zero for all $E\in \bR$.
\end{lemma}
\textbf{Proof:} This follows by symplectic reduction w.r.t.\ the free, proper 
Hamiltonian $\bR^d$--action on $P$, being the
symplectic lift of the diagonal $\bR^d$--action on $\widehat{M}$
(acting trivially on the momenta).
\hfill $\Box$\\[2mm]
Accordingly, we will work in the center of mass frame.

From a rough propagation estimate
it follows that for $j^+(x)=\infty$, the cluster function  $t\mapsto \cA(t)$, see \eqref{00V},
must change its value infinitely often.
We now show that for $n=4$ {\em there is a messenger cluster moving between two others}
infinitely often as $t\nearrow T(x)$.
\begin{enumerate}[1.]
\item 
{\em There is a nontrivial cluster:}
By Theorem \ref{thm:NAV:wandering}, for $x\in \operatorname{NAV}_{\!E,0}$ 
we have $\lim_{t\nearrow T(x)} \tilde{V}(t) = -\infty$.
Choose $t'_0\in(0, T(x))$ so that (with $\Xi_{\cC}^{(\delta)}$ from \eqref{Z4})
\[\tilde{V}(t) < \min \big\{ V(r) \mid r \in  (t'_0)^{1-\vep/2}\, \Xi_{\cC_{\wedge}}^{(\delta)} \big\}  
\qquad \big( t\in \big[t'_0, T(x)\big)\big).\] 
This is possible, since the minimum is increasing in the parameter $t'_0$.

Then by Definition \eqref{00V} of $\cA$, at least one atom of $\cA(t)$ must be nontrivial for all
$t\in [t_0, T(x))$, so that $|\cA(t)|\le3$.
\item 
{\em There is more than one cluster:}
On the other hand, since $j^+(x)=\infty$, 
$\lim_{t\nearrow T(x)} \|q(t,x)\|_\cM \, t^{\vep/2-1}=\infty$, too. 
As in the c.o.m.\ frame, the diameter of $\Xi_{\cC_{\vee}}^{(\delta)}$ is finite, we can
choose $t''_0\in(0, T(x))$ so that $\,t_0^{\vep/2-1}\not \in \Xi_{\cC_{\vee}}^{(\delta)} $ for $t\in [t''_0,T(x))$.
This implies that there are at least two clusters: $|\cA(t)|\ge2$.

For $t$ larger than $t_0:=\max(t'_0,t''_0)$ the number $|\cA(t)|$ of clusters is 2 or 3.
\item 
{\em There are infinitely many changes between two and three clusters:}
If for some time $\tau\in [t_0, T(x))$ one has $\cA(\tau^+)\neq \cA(\tau^-)$, then by construction of
the Graf partition \eqref{Z4}, 
$|\cA(\tau^+)| \neq |\cA(\tau^-)|$ and $\cA(\tau^+)$ is comparable with $\cA(\tau^-)$ in the sense of
Definition \ref{def:cluster}.

Breakup of a cluster followed by its recombination may happen. 

However, there must occur infinitely many 
changes between non-comparable set partitions, in the following sense. There exists a strictly increasing
sequence $(s_k)_{k\in\bN}$ of $s_k\in [s_0,T(x))$ with the properties $\cA(s_\ell^+)\neq \cA(s_\ell^-)$,
$\cA|_{(s_\ell,s_{\ell+1})}$ constant,
that $|\cA(s_\ell^+)|=3$ if $\ell$ is even and $|\cA(s_\ell^+)|=2$ if $\ell$ is odd, {\em and}
for every odd $\ell$ there is a minimal odd $\ell' > \ell$ with $\cA(s_{\ell'}^+)$ not comparable with 
$\cA(s_{\ell}^+)$. 

Otherwise $j^+(x)$ would be finite, since then the term 
$\frac{1}{2t^2} \| \tilde{Q}^E_\cA - \tilde{v}^E_\cA \|_\cM^2$ in 
\[\textstyle
\frac{d^2}{dt^2}\tilde{j}^E =  
-\frac{4}{t^2}\langle \tilde{Q}^E ,\frac{d}{dt}\tilde{q}^E - \tilde{Q}^E \rangle_\cM
+\frac{1}{2t^2} \| \tilde{Q}^E_\cA - \tilde{v}^E_\cA \|_\cM^2
-\frac{2}{t^2}   \langle \tilde{q}^E ,\nabla \tilde{V}^E \rangle\]
would be twice integrable, the two other terms having this property unconditionally 
(compare with $\frac{d}{dt}\tilde{j}^E$ in \eqref{djE}).
 \item
{\em Symbolic description:}
From this sequence we extract a subsequence $(t_k)_{k\in\bN}$, $t_k=s_{\ell(k)}$, maximal
with  respect to the properties that for $k$ even, 
\begin{enumerate}[$\bullet$]
\item 
$\cA$ is constant in the interval $(t_k,t_{k+1})$, and of size $|\cA(t_k^+)|=3$, 
\item 
$\ell(k\pm 1)=\ell(k)\pm 1$, $|\cA(t_{k}^-)|=|\cA(t_{k+1}^+)|=2$
and $\cA(t_{k}^-)$ and $\cA(t_{k+1}^+)$ are not comparable.
\end{enumerate}
Note that then for $k$ odd, $\cA$ need not  
to be constant in the interval $(t_k,t_{k+1})$, but
by maximality $\cA(t_k^+)= \cA(t_{k+1}^-)$, consisting of two clusters.

For $k$ even we consider the three-tuple $(\cC_{k-1},\cC_{k},\cC_{k+1})$ with $\cC_\ell:= \cA(t_\ell^+)$.
So $|\cC_{k}|=3$, $|\cC_{k\pm1}|=2$, $\cC_{k}$ is comparable with both $\cC_{k\pm1}$, but 
$\cC_{k-1}$ is not comparable with $\cC_{k+1}$. Thus we can uniquely enumerate the clusters
of $\cC_{k}$ in such a way that for $\cC_{k}=\{C_1,C_2,C_3\}$ we have 
\[\cC_{k-1}=\{C_1\cup C_2,C_3\} \qmbox{and} \cC_{k+1}=\{C_1, C_2 \cup C_3\}.\]
So $C_2$ is a cluster of particles whose center $\tilde{q}_{C_2}$
moves in the interval $(t_k,t_{k+1})$ from a neighborhood of  $\tilde{q}_{C_1}(t_k)$ to one
of  $\tilde{q}_{C_3}(t_{k+1})$.
It is called the {\em messenger}.
\item \label{combinatorics}
{\em Combinatorics:}
For $n=4$ the set ${\cal Q}$ of three three-tuples $(C_1,C_2,C_3)$ of clusters 
with $\{C_1,C_2,C_3\}\in \cP_3(N)$ (see \eqref{Pk}) \label{cal:T}
is of size $|{\cal Q}|={4\choose 2}3!=36$.
\end{enumerate}
Propagation estimates show (see Section \ref{sub:sec:prop:est}) that between their near-collisions, 
the three cluster centers move asymptotically as $t\nearrow T(x)$ on straight lines in~$\bR^d$. 

Thus we first consider in the following section \ref{sect3}
a nondeterministic kinematical model of three particles moving on straight lines
between their collisions. Then in Section \ref{sect4} we show that asymptotically 
the dynamics is indeed captured by that model.
Finally, in Section \ref{sub:sect:NAV:transitional} we define Poincar\'e surfaces which  
fulfill Assumptions {(A1)} and {(A2)} (Lemma \ref{lem:A1:A2}),  with
$\operatorname{NAV}_E\subseteq\operatorname{Trans}_E$.
%
\section{A nondeterministic kinematical model} \label{sect3}
%
Asymptotically, as $t\nearrow T(x)$, between their near-collisions the three 
cluster centers $\tilde{q}_{C_i}$ move on straight lines.
So to understand the kinematics, we set up a simple model, partly resembling the one of \cite{FKM}.

\begin{enumerate}[$\bullet$]
\item 
Three points with positions $\tilde{q}_1\le \tilde{q}_2\le \tilde{q}_3\in C^0 \big([0,T),\bR^d \big)$  
move with constant velocities $\tilde{v}_i=\frac{d}{dt} \tilde{q}_i$  
in configuration space $\bR^d$, until exactly two of them (say, $i$~and $j$) 
collide at collision times $t_k\in [0,T)$, that is, $\tilde{q}_i(t_k)=\tilde{q}_j(t_k)$.
\item 
The number of collisions happening in the future is infinite. Collisions are enumerated by 
$k\in \bN$, with $t_1:=0$, $t_{k+1}>t_k$ and $\lim_{k\to\infty}t_k=T$.\\ 
For $k$ even (odd), particle 2 collides with 3 (respectively 1) at time $t_k$.
\item 
Their masses may depend on time: 
\[\widetilde{m}_i:{\cal T}\to [m_{\min}, m_{\max}]\quad\mbox{with } 
{\cal T}:=[0,T)\setminus \cup_{k\in \bN}\{t_k\} \mbox{ and } m_{\min}>0.\] 
They are constant between collisions. The masses of the colliding particles 
may change during collision, but their sum is conserved:
\beq\widetilde{m}_i(t_k^+)+\widetilde{m}_j(t_k^+) = \widetilde{m}_i(t_k^-)+\widetilde{m}_j(t_k^-).
\Leq{geo:mass}
So the total mass $M:=\widetilde{m}_1+\widetilde{m}_2+\widetilde{m}_3$ is constant.
\item 
Denoting by $\tilde{p}_i:=\widetilde{m}_i\tilde{v}_i$ their momenta,  
during collision of particles $i$ and $j$, their external momentum is constant:
\beq
\tilde{p}_i(t_k^+)+\tilde{p}_j(t_k^+) = \tilde{p}_i(t_k^-)+\tilde{p}_j(t_k^-).
\Leq{geo:momentum} 
So by \eqref{geo:mass} at time $t_k$ their center of mass moves with constant velocity.
But the internal momentum 
$\frac {\tilde{m}_j\tilde{p}_i-\tilde{m}_i\tilde{p}_j}{\tilde{m}_i+\tilde{m}_j}$
may change arbitrarily, without conservation of total energy 
\beq
\textstyle \tilde{K}:= \sum_{i=1}^3 \frac{\|\tilde{p}_i\|^2}{2\widetilde{m}_i} : {\cal T}\to\bR. 
\Leq{def:K:tilde}
\item 
The initial conditions are chosen so that for $k=1$, that is, $t_1=0$ one has
$\tilde{J}'(t_k^+) > 0$ for the moment of inertia
\beq 
\textstyle 
\tilde{J} := \eh\sum_{i=1}^3 \widetilde{m}_i \|\tilde{q}_i\|^2:[0,T)\to \bR.
\Leq{ndm:initial}
\end{enumerate}
For application of the model to non-collision singularities, we would have $T\in (0,+\infty)$. For modeling
orbits having no asymptotic velocity, we take $T:=+\infty$.

The model indeed catches some properties of the 4-body orbits
whose asymptotic velocity does not exist. For example, particle 2 models the messenger cluster $C_2$
moving between $C_1$ and $C_3$. Of which of the four physical particles  these three clusters
consist, may change during near-collisions between three particles 
(so one would set $\widetilde{m}_i(t):=m_{C_i(t)}$).
Also, the cluster $D(t)$ consisting
of two particles serves as an infinite reservoir of kinetic cluster energy. 
\begin{prop}[Nondeterministic kinematical model] \quad\\[-6mm] 
\begin{enumerate}[1.]
\item 
The total angular momentum $\tilde{L}:=\sum_{i=1}^3 \tilde{L}_i$ with $\tilde{L}_i:=\tilde{q}_i\wedge \tilde{p}_i$
is conserved. In the center of mass frame, $\tilde{L}$ is zero or of rank two.
So the motion takes place on a line or a two-plane in configuration space $\bR^d$.

Considered as functions $\tilde{L}_i:{\cal T}\to\bR$, they are locally constant, and
\beq
|\tilde{L}_i|\le |\tilde{L}|\qquad (i=1,2,3).
\Leq{Li:bound}
\item
The moment of inertia $\tilde{J}$, see \eqref{ndm:initial}, is in $C^1([0,T),\bR)$, 
with $\frac{d^2}{dt^2}\tilde{J}(t)=2\tilde{K}(t)>0$ for $t\in {\cal T}$. Thus
\beq
\textstyle  \tilde{J}'(t_{k+1}^+)\ =\ \tilde{J}'(t_k^+)+2(t_{k+1}-t_k)\tilde{K}(t_k^+).
\Leq{F:prime:grows}
Moreover, for $\lambda:= \big(1+\frac{m_{\min}}{m_{\max}}\big)^{1/2}>1$ and some $J_0>0$
\beq
 \tilde{J}(t_k)\ge \lambda^k\,J_0\qquad (k\in \bN). 
\Leq{F:lower:bound}
\item
The speeds $\|\tilde{v}_i(t_k^+)\|$ are non-zero for all $k\in\bN$, and
\beq
\langle \tilde{p}_1(t_k^+),  \tilde{q}_1(t_k) \rangle > 0
\ \,  ,\ \, \langle \tilde{p}_2(t_k^+),  \tilde{q}_2(t_k) \rangle < 0
\ \mbox{ and }\
\langle \tilde{p}_3(t_k^+),  \tilde{q}_3(t_k) \rangle > 0.
\Leq{p:q}
\item 
The total kinetic energy $\tilde{K}(t_k^+)$ goes to infinity w.r.t.\ the collision index~$k$:
For $\mu:= \textstyle1+\big(\frac{m_{\min}}{m_{\max}} \big)^2 > 1$ and some $K_0>0$ 
depending on the initial conditions 
\beq
\tilde{K}(t_k^+)\ge \mu^k\, K_0 \qquad(k\in\bN).
\Leq{en:expo}
\item 
The directions $\hat{v}_i:=\tilde{v}_i/\|\tilde{v}_i\|$ align in the following sense: The limits 
\[w_i^{e}:= \lim_{k\to\infty} \hat{v}_i(t_{2k}^+)\in S^1 \qmbox{and}
w_i^{o}:= \lim_{k\to\infty} \hat{v}_i(t_{2k+1}^+)\in S^1\qquad (i=1,2,3)\] 
exist (index indicating {\em even} and {\em odd}), 
and $\,w_1^{e}=w_1^{o}=w_2^{e}=-w_2^{o}=w_3^{o}=w_3^{e}$.\\ 
The speed of convergence is exponential in $k$.
\end{enumerate}
\end{prop}
\begin{remark}[A tale of two scales]\quad\\ 
The exponential lower bounds \eqref{F:lower:bound} and \eqref{en:expo} depend on
collision index $k\in \bN$ instead of
time $t_k$. From \eqref{F:prime:grows} one infers that,
depending on the sequence $k\mapsto $ $t_{k+1}-t_k$ of time differences, the true behavior of these functions
can be super-exponential in $k$, and the growth of spatial extension is not coupled to the increase in speed. 
Celestial bodies on a line show the same behavior, see \cite{SX}.~$\Diamond$
\end{remark}
\textbf{Proof:}\\[-6mm]
\begin{enumerate}[1.]
\item 
$\tilde{L}$ is constant, since between collisions its derivative vanishes, and during collision of
particle $i$ and $j$ at time $t_k$ we have $\tilde{q}_i(t_k)=\tilde{q}_j(t_k)$ so that
\begin{align}
\label{ang:moment:model}
\tilde{q}_i(t_k)\wedge & \tilde{p}_i(t_k^-)+\tilde{q}_j(t_k)\wedge \tilde{p}_j(t_k^-) 
= \tilde{q}_i(t_k)\wedge \big(\tilde{p}_i(t_k^-)+\tilde{p}_j(t_k^-)\big)\\
&\!\!= \tilde{q}_i(t_k)\wedge \big(\tilde{p}_i(t_k^+)+\tilde{p}_j(t_k^+)\big)
= \tilde{q}_i(t_k)\wedge \tilde{p}_i(t_k^+)+\tilde{q}_j(t_k)\wedge \tilde{p}_j(t_k^+).\NN
\end{align}
For $k$ odd, at times $t_k<t_{k+1}<t_{k+2}$ particle two collides first with particle one, then
with three and again with one. In the inertial frame where $\tilde{q}_1(t)=0$ for $t\in (t_k,t_{k+2})$,
the momenta $\tilde{p}_2(t_k^+)$ and $\tilde{p}_3(t_k^+)$ span a line or plane. For having the collision
at time $t_{k+2}$, $\tilde{p}_2(t_{k+1}^+)$ must be linear dependent on 
$\tilde{p}_2(t_{k+1}^-)=\tilde{p}_2(t_k^+)$.
So also $\tilde{p}_3(t_{k+1}^+)$ must be in that line or plane, as is the center of mass. 
The case of $k$ even is similar. Thus ${\rm rank}(\tilde{L})\le2$.

Next we derive \eqref{Li:bound}. It is clear that the $\tilde{L}_i$ are locally constant, so that
$\tilde{L}_i(t_{k+1}^-) = \tilde{L}_i(t_{k}^+)$. As real-valued functions on $T^*E\cong \bR^2\times \bR^2$, 
the $L_i$ have the presentations $L_i = \langle{\mathbb J} q_i , p_i \rangle $, with 
${\mathbb J}:=\bsm 0 &-1\\ 1 & 0\esm$. 
We consider the case of $k\in\bN$ odd (the case of even $k$ being similar).
\begin{enumerate}
\item 
Solving for the momentum of the messenger particle two,
\beq \textstyle
\tilde{p}_2(t_k^+) = \widetilde{m}_2(t_k^+)\frac{\tilde{q}_3(t_k)-\tilde{q}_2(t_k)}{t_{k+1}-t_k} + 
\frac{\widetilde{m}_2(t_k^+)}{\widetilde{m}_3(t_k^+)}\, \tilde{p}_3(t_k^+).
\Leq{part:2:time}
We have $\tilde{q}_1(t_k) = \tilde{q}_2(t_k)$. Therefore in the center of mass frame 
$\tilde{q}_3(t_k)$ is a multiple of $\tilde{q}_2(t_k)$. So the first term
in \eqref{part:2:time} does not contribute to $\tilde{L}_2(t_{k}^+)$, and we get
\[ \textstyle
\tilde{L}_2(t_{k}^+) = \frac{\widetilde{m}_2(t_k^+)}{\widetilde{m}_3(t_k^+)} 
\langle {\mathbb J}\tilde{q}_2(t_k^+) , \tilde{p}_3(t_k^+) \rangle
= -  \frac{  \widetilde{m}_2(t_k^+) } {\widetilde{m}_1(t_k^+) + \widetilde{m}_2(t_k^+) } \tilde{L}_3(t_{k}^+)  .\]
\item 
Concerning particle one, we use that $\tilde{q}_1(t_k) = \tilde{q}_2(t_k)$ and
$\tilde{p}_1=-\tilde{p}_2-\tilde{p}_3$:
\begin{align}
\tilde{L}_1(t_{k}^+)& = \langle {\mathbb J}\tilde{q}_1(t_k^+) , \tilde{p}_1(t_k^+) \rangle
= - \tilde{L}_2(t_{k}^+)- \langle {\mathbb J}\tilde{q}_1(t_k^+) , \tilde{p}_3(t_k^+) \rangle\NN\\
&= \textstyle
\frac{  \widetilde{m}_2(t_k^+) +\widetilde{m}_3(t_k^+)  } {\widetilde{m}_1(t_k^+) + \widetilde{m}_2(t_k^+) } 
 \, \tilde{L}_3(t_{k}^+).
\end{align}
\item 
Thus  $\tilde{L} = \tilde{L}_1(t_{k}^+) +\tilde{L}_2(t_{k}^+)+ \tilde{L}_3(t_{k}^+)
= \frac{M}{M-\widetilde{m}_3(t_k^+)}  \tilde{L}_3(t_{k}^+)$, or
\beq
\textstyle
\tilde{L}_1(t_{k}^+) = \frac{M-\widetilde{m}_1(t_k^+)}{M}  \tilde{L}\mbox{ , }
  \tilde{L}_2(t_{k}^+) = - \frac{\widetilde{m}_2(t_k^+)}{M} \tilde{L}\mbox{ and }
  \tilde{L}_3(t_{k}^+) =  \frac{M-\widetilde{m}_3(t_k^+)}{M} \tilde{L}.
\Leq{3L:L}
\end{enumerate}
In particular, if the angular momentum $\tilde{L}(0)=0$ then for all $k$ the affine lines 
\[\ell_i(k):=\tilde{q}_i(t_k^+) + \operatorname{span}(\tilde{p}_i(t_k^+))\quad (i=1,2,3),\] 
on which the particles move in the time interval $[t_k,t_{k+1}]$, coincide.

\item
For non-collision times $t\in{\cal T}$ 
\beq
\textstyle\frac{d}{dt}\tilde{J}(t)=\sum_{i=1}^3   \langle\tilde{q}_i,\tilde{p}_i\rangle(t) \qmbox{and}
\frac{d^2}{dt^2}\tilde{J}(t)=2\tilde{K}(t). 
\Leq{F:derivatives}
By \eqref{geo:momentum}, for collision of particles $i$ and $j$ at time $t_k$
\beq
\langle \tilde{p}_i(t_k^+), \tilde{q}_i(t_k) \rangle + \langle \tilde{p}_j(t_k^+), \tilde{q}_j(t_k)\rangle
= \langle \tilde{p}_i(t_k^-), \tilde{q}_i(t_k)\rangle + \langle \tilde{p}_j(t_k^-), \tilde{q}_j(t_k)\rangle.
\Leq{F:prime:coll}
Thus we can uniquely extend its derivative so that $\tilde{J}\in C^1\big([0,T),\bR \big)$,
and its second derivative exists on ${\cal T}$ and is locally constant.

Equation \eqref{F:prime:grows} follows with \eqref{F:prime:coll} for all $k\in \bN$, 
since by \eqref{F:derivatives}
\[\textstyle \sum_{i=1}^3 \langle \tilde{p}_i(t_{k+1}^-), \tilde{q}_i(t_{k+1}) \rangle 
=\textstyle \sum_{i=1}^3 \langle \tilde{p}_i(t_k^+), \tilde{q}_i(t_k) \rangle +2(t_{k+1}-t_k)\tilde{K}(t_k^+),\]
and $\tilde{K}(t_k^+)>0$, since otherwise there would not be a further collision.

To prove \eqref{F:lower:bound}, we first analyze the special case of motion on a line.
Then both \eqref{geo:momentum},  the motion between collisions
\[\textstyle \tilde{q}_i(t_{k+1})
= \tilde{q}_i(t_k) + \frac{\Delta t_k}{\widetilde{m}_i(t_k^+)}\tilde{p}_i(t_k^+) \qquad (k\in \bN, i=1,2,3)\]
and the collision conditions 
\[\tilde{p}_2(t_{k}^+) = \frac{\widetilde{m}_2(t_{k}^+)}{\Delta t_k}\times  \left\{
\begin{array}{cc}
(\tilde{q}_3(t_{k+1})-\tilde{q}_1(t_k)) &,\, k\mbox{ odd}\\
(\tilde{q}_1(t_{k+1})-\tilde{q}_3(t_k)) &,\, k\mbox{ even}
\end{array}\right.\]
at times $t_k$ are linear in the $\tilde{q}_i(t_k)$ 
and $\tilde{p}_i(t_k^+)$, the coefficients being rational functions in the time differences 
$\Delta t_k:=t_{k+1}-t_k>0$
and the masses $\widetilde{m}_i(t_k^+)$.
For $\Delta \tilde{q} :=\tilde{q}_3 -\tilde{q}_1$ we obtain for two subsequent collisions
a somewhat large set of linear equations in the positions and momenta at times $t_k$ and 
$t_{k+2}$.
As we know the sign of its coefficients, 
assuming \eqref{p:q} (whose proof will be not use estimate \eqref{F:lower:bound}),
without writing the equations explicitly, we have
\[\Delta \tilde{q}(t_{k+2}) \ge \left\{
\begin{array}{cc}
\big(1+\frac{\widetilde{m}_2(t_{k+1}^+)}{\widetilde{m}_1(t_{k+1}^+)}\big)\Delta \tilde{q}(t_k) &,\, k\mbox{ odd}\\
\big(1+\frac{\widetilde{m}_2(t_{k+1}^+)}{\widetilde{m}_3(t_{k+1}^+)}\big)\Delta \tilde{q}(t_k) &,\, k\mbox{ even}
\end{array}\right. .\]
Since $\Delta \tilde{q}(t_{1})>0$, this proves estimate \eqref{F:lower:bound}, with $\lambda =
 \big(1+\frac{m_{\min}}{m_{\max}}\big)^{1/2}$.

For the generic case of non-vanishing angular momentum $\tilde{L}$,
we already know that the motion in the plane is asymptotic to a one on a line. 
So by a continuity argument we get the result also in the general case.
\item
That the speeds $\|\tilde{v}_i(t_k^+)\|$ are non-zero follows from the inequalities \eqref{p:q}.
We show that these are valid for all $k\in\bN$, assuming w.l.o.g.\ that $k$ is odd.
\begin{enumerate}
\item 
As $ \tilde{q}_1(t_k)= \tilde{q}_2(t_k)$, in the center of mass system
\beq
0 < \tilde{J}'(t_k) = \textstyle \sum_{i=1}^3 \langle \tilde{p}_i(t_k^+), \tilde{q}_i(t_k)\rangle =
\langle \tilde{p}_3(t_k^+), \tilde{q}_3(t_k)\rangle \frac{M} {\widetilde{m}_3(t_k^+)},
\Leq{p3q3:g0}
showing the third inequality in \eqref{p:q}. Additionally we get for $i=1,2$
\beq
 \textstyle \tilde{q}_i(t_k) = -\frac{\widetilde{m}_3(t_k^+)}{M-\widetilde{m}_3(t_k^+)} \tilde{q}_3(t_k^+)
\qmbox{and thus} 
\langle \tilde{p}_3(t_k^+), \tilde{q}_i(t_k)\rangle =- \frac{M\,\tilde{J}'(t_k)}{M-\widetilde{m}_3(t_k^+)} .
\Leq{scalar:products}
\item 
We now take the next collision at time $t_{k+1}$ and between particles 2 and 3 into account.
By \eqref{part:2:time} and both identities in \eqref{scalar:products}
\[\textstyle
\langle \tilde{p}_2(t_k^+), \tilde{q}_2(t_k)\rangle = 
- \frac{\widetilde{m}_2 M}{M-\widetilde{m}_3}
\left(\frac{\widetilde{m}_3}{M-\widetilde{m}_3} \frac{\|\tilde{q}_3\|^2}{t_{k+1}-t_k} 
+\frac{\tilde{J}'}{\widetilde{m}_3}\right)(t_k^+)<0.\]
\item 
To prove the first inequality in \eqref{p:q}, we note that by the above
\beq
\langle \tilde{p}_1(t_k^+),  \tilde{q}_1(t_k) \rangle = 
- \langle \tilde{p}_2(t_k^+),  \tilde{q}_2(t_k) \rangle
- \langle \tilde{p}_3(t_k^+),  \tilde{q}_2(t_k) \rangle >0.
\Leq{p1q1:g0}
\end{enumerate}
\item
With $\hat{q}_i:=\tilde{q}_i/\|\tilde{q}_i\|$, at collision times $t_k$ we estimate from below the terms 
\beq
\tilde{K}_\parallel:{\cal T}\to\bR\qmbox{,}
\tilde{K}_\parallel:= 
\frac{\langle \tilde{p}_1,\hat{q}_1\rangle^2+\langle \tilde{p}_3,\hat{q}_3\rangle^2}{2M} 
\Leq{K:parallel}
in the kinetic energies, corresponding to the momentum components
that are then parallel to the line through the positions of all three particles.

We know from \eqref{p:q} that $\tilde{K}_\parallel(t_1^+)>0$. Also
$\tilde{K}_\parallel<\tilde{K}$, since, unlike in \eqref{def:K:tilde}, the mass $M$ appears in \eqref{K:parallel}, 
and $\widetilde{m}_i<M$.
We prove that
\beq
\tilde{K}_\parallel(t_{k+1}^+)\ge \mu\, \tilde{K}_\parallel(t_k^+)\qmbox{with}
\mu = \textstyle 1+\big(\frac{m_{\min}}{m_{\max}}\big)^2 > 1,
\Leq{K:mu:K}
assuming w.l.o.g.\ that $k$ is odd. Then \eqref{en:expo} follows with $K_0:=\tilde{K}_\parallel(t_1^+)/\mu$.
\begin{enumerate}
\item 
The term 
\[\textstyle C_1(k):=
\frac{\langle \tilde{p}_1(t_{k+1}^+),\hat{q}_1(t_{k+1})\rangle^2}{2M}
- \frac{\langle \tilde{p}_1(t_{k}^+),\hat{q}_1(t_{k})\rangle^2}{2M}\]
in $\tilde{K}_\parallel(t_{k+1}^+) - \tilde{K}_\parallel(t_k^+)$ is positive, since 
$\tilde{p}_1(t_{k+1}^+) = \tilde{p}_1(t_{k}^+)$ and
\[\textstyle 
\tilde{q}_1(t_{k+1})=\tilde{q}_1(t_{k})+\frac{\tilde{p}_1(t_{k}^+)}{\widetilde{m}_1(t_k^+)}(t_{k+1}-t_k)
\ \mbox{ with } t_{k+1}-t_k>0, \]
so that $\hat{q}_1(t_{k+1})$ is more aligned to $\tilde{p}_1(t_{k}^+)$ than
$\hat{q}_1(t_k)$.
\item 
The second term in  $\tilde{K}_\parallel(t_{k+1}^+) - \tilde{K}_\parallel(t_k^+)$,
\[\textstyle C_3(k):=
\frac{\langle \tilde{p}_3(t_{k+1}^+),\hat{q}_3(t_{k+1})\rangle^2}{2M}
- \frac{\langle \tilde{p}_3(t_{k}^+),\hat{q}_3(t_{k})\rangle^2}{2M},\]
is more complicated and will be estimated more precisely. We have
\beq \tilde{p}_3(t_{k+1}^+) = \tilde{p}_3(t_{k+1}^-)+ \big(\tilde{p}_2(t_{k+1}^-) - \tilde{p}_2(t_{k+1}^+)\big)
\Leq{total:momentum:conserved}
and $\tilde{p}_i(t_{k+1}^-)= \tilde{p}_i(t_k^+)$. Thus 
\beq 
C_3(k) = \frac{\big(D_{3,I}^a(k) + D_{3,I\!I}(k) + D_{3,I\!I\!I}(k) \big)^2 - \big(D_{3,I}^b(k) \big)^2}{2M}
\Leq{C:CCC}
with
\begin{align*}
\textstyle 
&D_{3,I}^a(k)\ \, := 
 \langle \tilde{p}_3(t_k^+),\hat{q}_3(t_{k+1}) \rangle\mbox{ , }
 D_{3,I}^b(k) := 
 \langle \tilde{p}_3(t_k^+), \hat{q}_3(t_{k}) \rangle\\
&D_{3,I\!I}(k)\; := \langle \tilde{p}_2(t_k^+),\hat{q}_3(t_{k+1})  \rangle\qmbox{and}\\
&D_{3,I\!I\!I}(k):= -\langle \tilde{p}_2(t_{k+1}^+),\hat{q}_3(t_{k+1})  \rangle
.
\end{align*}
\begin{enumerate}[I.]
\item 
By \eqref{p3q3:g0} and the same argument as in (a), one has 
\[D_{3,I}^a(k)\ge D_{3,I}^b(k)>0.\]
\item 
By solving the identity
$\frac{\tilde{p}_2(t_k^+)}{\widetilde{m}_2(t_k^+)}(t_{k+1}-t_k)=\tilde{q}_2(t_{k+1})-\tilde{q}_2(t_k)
=\tilde{q}_3(t_{k+1})-\tilde{q}_1(t_k)=\tilde{q}_3(t_k)-\tilde{q}_1(t_k)+ \frac{\tilde{p}_3(t_k^+)}{\widetilde{m}_3(t_k^+)}(t_{k+1}-t_k)$ for $\tilde{p}_2(t_k^+)$,
\begin{align*}
D_{3,I\!I}(k) &=\textstyle \widetilde{m}_2(t_k^+)
\big\langle \frac{\tilde{q}_3(t_k)-\tilde{q}_1(t_k)}{t_{k+1}-t_k} 
  + \frac{\tilde{p}_3(t_k^+)}{\widetilde{m}_3(t_k^+)},\hat{q}_3(t_{k+1})  \big\rangle \\
  &>\textstyle\frac{\widetilde{m}_2(t_k^+)}{\widetilde{m}_3(t_k^+)}
\langle  \tilde{p}_3(t_k^+),\hat{q}_3(t_{k+1})  \rangle 
>\textstyle\frac{\widetilde{m}_2(t_k^+)}{\widetilde{m}_3(t_k^+)}
\langle  \tilde{p}_3(t_k^+),\hat{q}_3(t_k)  \rangle \; > 0.
\end{align*}
\item 
Finally we use \eqref{p1q1:g0} to show that
\begin{align*}
D_{3,I\!I\!I}(k)&=\textstyle \widetilde{m}_2(t_{k+1}^+)
\big\langle \frac{\tilde{q}_3(t_{k+1})-\tilde{q}_1(t_{k+1})}{t_{k+2}-t_{k+1}} 
  - \frac{\tilde{p}_1(t_{k+1}^+)}{\widetilde{m}_1(t_{k+1}^+)},\hat{q}_3(t_{k+1})  \big\rangle \\
  & > \textstyle - \frac{\widetilde{m}_2(t_{k+1}^+)}{\widetilde{m}_1(t_{k+1}^+)} 
  \big\langle \tilde{p}_1(t_{k+1}^+),\hat{q}_3(t_{k+1}) \big\rangle\\
  & = \textstyle \frac{\widetilde{m}_2(t_{k+1}^+)}{\widetilde{m}_1(t_{k+1}^+)} 
  \big\langle \tilde{p}_1(t_k^+),\hat{q}_1(t_{k+1}) \big\rangle\\
  & \ge \textstyle \frac{\widetilde{m}_2(t_{k+1}^+)}{\widetilde{m}_1(t_{k+1}^+)} 
  \big\langle \tilde{p}_1(t_k^+),\hat{q}_1(t_k) \big\rangle\; > 0.
\end{align*}
\end{enumerate}
Thus 
\eqref{C:CCC} is estimated from below, with $\tilde{K}_\parallel$ from \eqref{K:parallel}, by
\begin{align*}
C_3(k) &\textstyle  > \frac{(D_{3,I\!I}(k) )^2+ (D_{3,I\!I\!I}(k) )^2}{2M}\\
&> \textstyle\frac{\Big(\frac{\widetilde{m}_2(t_k^+)}{\widetilde{m}_3(t_k^+)}
\langle  \tilde{p}_3(t_k^+),\hat{q}_3(t_k)  \rangle\Big)^2
+ \Big(\frac{\widetilde{m}_2(t_{k+1}^+)}{\widetilde{m}_1(t_{k+1}^+)} 
  \big\langle \tilde{p}_1(t_k^+),\hat{q}_1(t_k) \big\rangle\Big)^2}{2M}\\
 & \ge \textstyle \big(\frac{m_{\min}}{m_{\max}}\big)^2\tilde{K}_\parallel(t_k^+).
\end{align*}
\end{enumerate}
Together, (a) and (b) prove \eqref{K:mu:K}, and thus \eqref{en:expo} follows.
\item
The alignment of the directions $\hat{v}_i$ follows from the above statements, using the bounds
on the angular momenta:

There is a unique closed circular segment 
$\operatorname{seg}(\hat{v}_2(t_{k+1}^-), \hat{v}_3(t_{k+1}^-))\subseteq S^1$ of length
in $[0,\pi)$ between $\hat{v}_2(t_{k+1}^-)$ and $\hat{v}_3(t_{k+1}^-)$. 
Still the oriented lines $\ell_1(k)$ and $\ell_3(k)$ may be antiparallel, 
so that $\hat{v}_1(t_{k+1}^-) =- \hat{v}_3(t_{k+1}^-)$.
One concludes from $t_{k+1}\in (t_k,+\infty)$ that for $k$ odd
\beq
-\hat{v}_2(t_{k+1}^+)\in\left\{
\begin{array}{cc}
\operatorname{seg}\big(\!-\!\hat{v}_1(t_{k+1}), \hat{v}_2(t_{k+1}^-)\big) &,\, k\mbox{ odd}\\
\operatorname{seg}\big(\!-\!\hat{v}_3(t_{k+1}), \hat{v}_2(t_{k+1}^-)\big) &,\, k\mbox{ even}
\end{array}
\right. ,
\Leq{ptwo:in:segments}
see Figure \ref{alignment}. The subsets 
\beq
\hspace*{-2.5mm}
S_k:=\left\{
\begin{array}{cc}
\!\!\operatorname{seg}\big(\!-\!\hat{v}_1(t_{k+1}), \hat{v}_2(t_{k+1}^-)\big) \cup\, 
\operatorname{seg}\big(\hat{v}_2(t_{k+1}^-), \hat{v}_3(t_{k+1}^-)\big) \!\!\!&,\, k\mbox{ odd}\\
\!\!\operatorname{seg}\big(\!-\!\hat{v}_3(t_{k+1}), \hat{v}_2(t_{k+1}^-)\big) \cup\, 
\operatorname{seg}\big(\hat{v}_2(t_{k+1}^-), \hat{v}_1(t_{k+1}^-)\big) \!\!\!&,\, k\mbox{ even}
\end{array}
\right.
\Leq{def:Sk} 
of $S^1$ are segments, 
since both segments on the r.h.s.\ of \eqref{def:Sk} contain the end point $\hat{v}_2(t_{k+1}^-)$.

From the preservation of total momentum at time $t_{k+1}$, that is \eqref{total:momentum:conserved},
it follows that 
\beq
S_{k+2}\subseteq S_k \qquad(k\in \bN).
\Leq{S:inclusion}

Moreover, we show now that the lengths of the segments $S_k$ converge to zero exponentially in $k$:
\begin{enumerate}[$\bullet$]
\item 
At the times $t_k$ of collision, in the center of mass frame, all three particles are on a line through
the origin, and from \eqref{F:lower:bound}, exponentially in $k$
\[\lim_{k\to\infty} \|\tilde{q}_i(t_k)\| = \infty \qquad (i=1,2,3).\]
\item 
By \eqref{en:expo} the last statement is also true for the speeds: Exponentially in $k$
\[\lim_{k\to\infty} \|\tilde{v}_i^\pm(t_k)\|=\infty \qquad (i=1,2,3).\]
\item 
On the other hand, the angular momenta are uniformly bounded by \eqref{Li:bound}.
\end{enumerate}
These three statements are only compatible, if the directions $\hat{v}_i(t_k^\pm)\in S^1$  ($i=1,2,3$)
become (anti-) parallel as $k\to \infty$. Together with \eqref{S:inclusion}, we conclude that, as $k\to \infty$, 
there is a limiting line they span.
\hfill $\Box$
\end{enumerate}
%

%
\begin{figure}[htbp]
\begin{center}
\includegraphics[height=35mm]{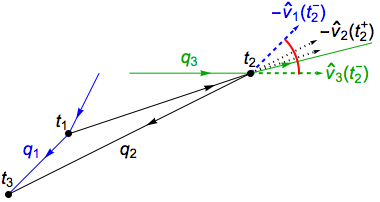}
\hfill
\includegraphics[height=35mm]{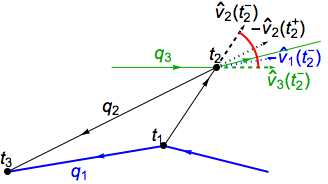}
\caption{Alignment of velocities, here of particle 2 and 3 at time $t_{k+1}\equiv t_2$. 
The segment $S_k$ containing $\hat{v}_3(t_2^+)$ is shown in red. Note the different orientations in the
two pictures, and that the third possibility
$S_k=\operatorname{seg}\big(\!-\!\hat{v}_1(t_{k+1}^-), \hat{v}_2(t_{k+1}^-)\big)$
is also realizable.}
\label{alignment}
\end{center}
\end{figure}
%
%
\section{Dynamics} \label{sect4}
%
The aim of this section is to provide the dynamical estimates
that will be used to show, in Section \ref{sub:sect:NAV:transitional}, that for every total energy $E\in\bR$, 
with the ${\frak L}$-dependent family of Poincar\'e sections $\cH_{m,E,{\frak L}}  \subseteq \Sigma_{E,0}$
defined in \eqref{the:surfaces} and \eqref{HmC}, and
\beq
\operatorname{Trans}_{E,0,{\frak L}}:=
\{x\in \Sigma_{E,0} \mid 
\exists\, m_0\in\bN\;\forall \,m\ge m_0: \cO^+(x)\cap \overline{\cH}_{m,E,{\frak L}}\neq\emptyset \}\, , 
\Leq{trans:L}
\[\textstyle
\operatorname{NAV}_{E,0} \, \subseteq \,\bigcup_{{\frak L}>0} \operatorname{Trans}_{E,0,{\frak L}} \, .\]
In words: in the center of mass frame, for every initial condition $x$ of energy $E$
whose asymptotic velocity limit \eqref{asymptotic:velocity} does not exist,
there is an angular momentum parameter ${\frak L}$, such that its forward orbit intersects almost all 
Poincar\'e surfaces $\cH_{m,E,{\frak L}} $ of that parameter.

Considering the definition \eqref{HmC} of the Poincar\'e surfaces, the main task will be to control the
evolution of the cluster angular momenta.  Total angular momentum $L$, see 
\eqref{total:ang:momentum}, is a constant of the motion. As shown in Lemma \ref{lem:cluster:ang:mom}
(based on estimates of the current section)
$\|L(x)\|$ sets a scale for the angular momenta of the clusters, 
and thus for the parameter ${\frak L}$ for which $x\in \operatorname{Trans}_{E,0,{\frak L}}$.

It has been proven in Section \ref{sub:sect:B2} that for initial conditions $x\in\operatorname{NAV}$
and times $t\in [t_0,T(x))$ the cluster size $|\cA(t)|$ is two or three, 
and that the time interval is partitioned into subintervals $(t_k,t_{k+1})$ where 
for $k$ even, $\cA$ is constant and of size $|\cA(t_k^+)|=3$, whereas $|\cA(t_{k}^-)|=|\cA(t_{k+1}^+)|=2$
and $\cA(t_{k}^-)$ and $\cA(t_{k+1}^+)$ are not comparable.
Accordingly, we have to tackle two dynamical problems:
\begin{enumerate}[$\bullet$]
\item 
For $k$ odd, the internal motion of the cluster $C_1\cup C_2$, perturbed by cluster~$C_3$.
Then the relative cluster angular momentum \eqref{rel:cluster:ang}
of the pair $C_1$, $C_2$ of clusters is small, see Section \ref{sub:sec:two:clusters}.
\item 
For $k$ even and $\cA(t)=\{C_1,C_2,C_3\}$, the flight of the messenger cluster $C_2$
from $C_1$ to $C_3$. The cluster centers will asymptotically move on straight lines, 
see Section \ref{sub:sec:prop:est}.
\end{enumerate}
%
\subsection{Preparatory estimates} \label{sub:sec:prop:prep}
%

To obtain the dynamical estimates, we use: 
\begin{enumerate}[1.]
\item 
\label{dyn:assertion:1}
By the von Zeipel Theorem \ref{thm:von:Zeipel},
the limit $j^+(x):=\lim_{t\nearrow T^+(x)} j_x(t)$ exists and equals $+\infty$.
Additionally, inspection of \eqref{djE} shows that for all $t\in (t_k,t_{k+1})$ the time derivative  
$\frac{d}{dt}\tilde{J}^E_\cA$  of the external moment of inertia is positive and goes to infinity. 
For $k$ odd this means that both
$C_1\cup C_2$ and $C_3$ move away from the origin, with opposite momenta
(in the center of mass frame).
\item 
\label{dyn:assertion:2}
When $t$ is a 'messenger time', {\em i.e.}, $t\in \operatorname{mess}$ with
\[\operatorname{mess}:= \cup_{k\in 2\bN}(t_k,t_{k+1}), \]
we have $|\cA(t)|=3$ clusters, 
and the messenger cluster $C_2$, see above.
This implies for the unique nontrivial cluster $D(t)\in \cA(t)$, $t\in\operatorname{mess}$,  that 
\end{enumerate}
\begin{lemma}[The nontrivial cluster] For $x\in \operatorname{NAV}$ and $E:= H(x)$,
\label{lem:nontrivial:cluster}%
\begin{enumerate}[(a)]
\item 
The external potential energy 
$\ \tilde{V}^E_{\cA(t)}(t)  = o(1) \quad (t\nearrow T(x))$\,,\\ 
and $\tilde{V}^I_{\cA(t)}(t)=\tilde{V}^I_{D(t)}(t)$ for $t\in\operatorname{mess}$.
\item 
The external cluster energy has the limit
$\ \lim_{t\nearrow T(x)}\tilde{H}^E_{\cA(t)}(t) = +\infty$. 
\item 
The internal energy of the unique nontrivial cluster $D(t)$ has the limit
\[\textstyle
\lim_{\substack {t\nearrow T(x)\\ t\in\operatorname{mess}}}\ \tilde{H}^I_{D(t)}(t) = -\infty.\] 
\item 
The internal angular momentum (see \eqref{rel:ang:mom}) of $D(t)$ has the limit 
\[\textstyle
\lim_{\substack {t\nearrow T(x)\\ t\in\operatorname{mess}}}\ \tilde{L}^I_{D(t)}=0.\]  
\end{enumerate}
\end{lemma}
{\bf Proof:}
\begin{enumerate}[(a)]
\item 
The first property is valid for {\em all} initial conditions $x\in P$, 
as in $\cA(t)$ (see \eqref{00V}) the different cluster centers have distances 
bounded below by $\delta^2 t^{1-\vep/2}$.
The identity $\tilde{V}^I_{\cA(t)}(t)=\tilde{V}^I_{D(t)}(t)$ 
follows from the Definition \eqref{V:IE:C} of $V^I_\cC$, since
both clusters in $\cA(t)\setminus\{C(t)\}$ are trivial.
\item 
By Theorem \ref{thm:NAV:wandering}.1, $\lim_{t\nearrow T(x)}\|\frac{d}{dt}q_{\cA(t)}^E(t) \|_\cM=\infty$
for the cluster-external speed. So the total kinetic cluster energy diverges:
$\lim\limits_{t\nearrow T(x)}\tilde{K}^E_{\cA(t)}(t)=+\infty$.
Since $\tilde{H}^E_{\cA}=\tilde{K}^E_{\cA}+\tilde{V}^E_{\cA}$ (see \eqref{H:IE:C}), 
(b) follows from Part (a).
\item 
follows from Part (b) by energy conservation, 
as $E= \tilde{H} = \tilde{H}^E_\cA + \tilde{H}^I_\cA$, and for 
$t \in \operatorname{mess}$, $\tilde{H}^I_{\cA(t)}(t) = \tilde{H}^I_{D(t)}(t)$.
\item 
As $H^I_D = K^I_D + V^I_D$, with 
$K^I_D(p^I_D,q^I_D)\ge \frac{\|L^I_D\|^2}{2m^I_D\|q^I_D\|^2}$ and 
$V^I_D(q^I_D)\ge -\frac{I}{\|q^I_D\|^\alpha}$ with $\alpha<2$, Statement (c)
is only possible if Assertion (d) is true. 
\hfill $\Box$
\end{enumerate}

For the three-tuple $(C_1,C_2,C_3)\in{\cal T}$ with $\cA(t)=\{C_1,C_2,C_3\}$, 
attributed to the interval $(t_k,t_{k+1})$ (see page \pageref{cal:T}), 
the messenger cluster $C_2$ carries a positive proportion of
the total kinetic energy $\tilde{K}^E_{\cA} = \sum_{i=1}^3 \tilde{K}^E_{C_i}$:
\begin{lemma} [Kinetic energy of the messenger cluster] \label{lem:part:energy}
\[ \textstyle \tilde{K}^E_{C_2(t)}(t) \ge \frac{m_{\min}}{7\,m_{\max}} \tilde{K}^E_{\cA(t)}(t)
\qquad(t\in\operatorname{mess}).\]
\end{lemma}
\textbf{Proof:}
This follows for {\em a certain\,} $t\in (t_k,t_{k+1})$, since the speed $\|v_{C_2}\|$ of the messenger particle
must be at least the one of the cluster $C_3$ it is to reach, and
\begin{align*}
2\tilde{K}^E_{\cA}&= 
m_{C_1} \|\tilde{v}_{C_2}+\tilde{v}_{C_3}\|^2+m_{C_2} \|\tilde{v}_{C_2}\|^2+m_{C_3}\|\tilde{v}_{C_3}\|^2\\
&\le m_{C_1} \|2\tilde{v}_{C_2}\|^2\!+m_{C_2} \|\tilde{v}_{C_2}\|^2\!+m_{C_3}\|\tilde{v}_{C_2}\|^2 
\textstyle\le 6\,m_{\max}\|\tilde{v}_{C_2}\|^2\le \frac{6m_{\max}}{m_{\min}}2\tilde{K}^E_{C_2}.
\end{align*}
It will follow for {\em all\,} $t\in (t_k,t_{k+1})$ (with factor 7 instead of 6)
from the propagation estimates of Sect.\ \ref{sub:sec:prop:est}, since the velocities $\tilde{v}_{C_2}$
are nearly constant in $(t_k,t_{k+1})$.~$\Box$\\[2mm]
Thus, when at time $t_k$ cluster $C_2$ separates from cluster $C_1$, their {\em relative}
kinetic energy is not negligible, compared to the internal energy of the cluster $D$.
This fact will be used in Lemma \ref{lem:small:var:a:m}, concerning the motion of the messenger cluster.
Although the scheme $\ \ldots\ar\cA(t_k)\ar \cA(t_{k+1})\ar\ldots$, 
depicted in \eqref{diagram}, 
\beq
\hspace*{-2.2mm}\{C'_1,C'_2\cup C_3\}\!\rightarrow \! \{C'_1,C'_2,C_3\} \! \rightarrow 
\stackrel {{\color{red}\rm Section\ } \ref{sub:sec:two:clusters}}{\{C_{1,2},C_3\}}  \rightarrow 
\stackrel{{\color{red}\rm Section\ } \ref{sub:sec:prop:est}}{ \{C_1,C_2,C_3\} } \rightarrow \!\{C_1,C_2\cup C_3\}
\Leq{diagram}
repeats after two steps, when showing in Section \ref{sub:sec:two:clusters} 
that the relative angular momentum of $C_1$ and $C_2$ must be small we use that in both time
directions the messenger cluster $C_2$ resp.\ $C'_2$ returns from 
$C_{1,2}:=C_1\cup C_2=C'_1\cup C'_2$ to $C_3$.

As Sections  \ref{sub:sec:two:clusters} and \ref{sub:sec:prop:est}  consider the situation where 
the cluster centers of $C_2$ and $C_1$ (or $C_3$) are initially very close, we now introduce the 
quantities that are adapted to that situation.

For two mutually disjoint clusters $C,D\subseteq N$, the following quantities neither depend on
the reference frame nor on the ordering of $C$ and $D$:
total mass $m_{C\cup D} =m_C+m_D$, 
reduced mass $m_{C,D}:=\frac{m_C m_D}{m_C+m_D}$,
and the {\em relative} phase space functions 
\begin{align}
L_{C,D} &:= \eh (q_{C}-q_{D}) \wedge (p_{C}-p_{D}) 
\label{rel:cluster:ang}\\
K_{C,D}&:=  \textstyle \eh m_{C,D}\|v_C-v_D\|^2
\label{K:C2:Ci}\\
H_{C,D} &:= \textstyle K_{C,D} +V_{C,D} \qmbox{with} V_{C,D} := \sum_{i\in C, j \in D} V_{i,j}
\label{Ham:C:D}\\
\textstyle
J_{C, D}(q) &:=\eh(m_{C}\|q_{C}\|^2 + m_{D}\|q_{D}\|^2 ) .
\label{J:C:D}
\end{align}
%
\subsection{Near-collision of two clusters} \label{sub:sec:two:clusters}
%
We now consider the near-collision of the clusters $C_1$ and $C_2$
that took place before time $s_1$, see \eqref{diagram}.

For simplicity of presentation, with the involution ${\cal R}:(p,q)\mapsto (-p,q)$ 
we reverse the direction of time. 
That does not change the various kinetic energies and reverses the sign of the angular momenta.
W.l.o.g.\ we set $s_1:=0$.

If initially (at time $0$) the cluster $D\in \{C_1,C_2,C_3\}$ that is nontrivial, equals $C_1$ or $C_2$,
a close encounter of the three bodies $C_{1,2} = C_1 \cup C_2$ follows.\\ 
Otherwise $D=C_3$, and the close encounter involves only two bodies.

Since the initial condition is $x\in \operatorname{NAV}_E$, we can assume that 
\begin{enumerate}[1.]
\item 
$|E| \ll -\tilde{H}^I_D(s_1)$, as $\tilde{H}^I_D(t)\searrow -\infty$ by Lemma \ref{lem:nontrivial:cluster} (c).
\item 
The relative kinetic energy $K_{C_1,C_2}$
of the two clusters (see \eqref{K:C2:Ci}) is large: 
\[\textstyle
\tilde{K}_{C_1,C_2}(0)\ge \frac{m_{\min}}{7\,m_{\max}}\, ( - \tilde{H}^I_D(0)).\]
That assumption is justified, as by Lemma \ref{lem:part:energy} on the partition of kinetic energy
$ \textstyle \tilde{K}^E_{C_2} \ge \frac{m_{\min}}{7\,m_{\max}} \tilde{K}^E_{\cA}$.
\item 
Repeating the notation from \eqref{diagram}, after having reversed the direction of time,
\begin{enumerate}
\item 
in the past $C_2$ had a near-collision with $C_3$
\item 
after the near-collision between $C_1$ and $C_2$,  $C_{1,2}$  decomposes into
$C'_1$ and $C'_2$ ($C_{1,2}=C'_1 \cup C'_2$),
\item 
finally  $C'_2$ has a near-collision with $C'_3=C_3$.
\end{enumerate}
\end{enumerate}
As by the propagation estimates between near-collision the cluster centers have nearly constant
velocities, this means that the direction of $p_{C'_2}$ must be nearly opposite to the one of 
$p_{C_2}$. We show that this can only be the case if initially the norm of the 
relative angular momentum $L_{C_1,C_2}$, see \eqref{rel:cluster:ang},
of $C_1$ and $C_2$ has small modulus. 
The basic estimate is Statement \ref{direction} of Lemma \ref{lem:three:particles:no:int}. \\
First, in Lemma \ref{lem:three:particles:no:int}, we consider the case of {\em no interaction} between
the two clusters $C_{1,2}$ and~$C_3$:
$V^E_{\{C_{1,2},C_3\}}=0$,  see \eqref{V:IE:C} and Assumption \ref{lem:ass:3} below.
Then the centers of mass $q_{C_{1,2}}$ and $q_{C_3}$ move with constant velocities.

So we change from the center of mass frame of $C_N$ to the one of $C_{1,2}$, {\em i.e.}
\beq
\textstyle\sum_{i\in C_{1,2}} p_i=0 \qmbox{,} \sum_{i\in C_{1,2}}m_iq_i  = 0 \,,
\Leq{com:c12}
{\em without modifying the notation}. 
By Assumption \ref{lem:ass:3}, $H^I_{C_{1,2}}$ is a constant of the motion, and its Hamiltonian equations 
describe the relative motion of the single particles $C_1$ and $C_2$. 
If $D\neq C_3$,
\beq \textstyle
H^I_{C_{1,2}}=H_{C_1,C_2}+H^I_{D}
\qmbox{for}
H_{C_1,C_2} = K_{C_1,C_2} +V_{C_1,C_2},
\Leq{Ham:G}
see \eqref{Ham:C:D}.
For $D=C_3$, $H^I_{C_{1,2}}=H_{C_1,C_2}$, and $H^I_{D}$ is a constant of the motion.
The total angular momentum $L_{C_{1,2}} = \sum_{i=1}^3 q_i\wedge p_i $
of $C_{1,2}$ is a constant of motion, too.

When we treat the case $D=C_1$ (or equivalently $D=C_2$), 
by permuting indices, if necessary, we assume $C_1=D=\{1,2\}$, $C_2=\{3\}$ and $C_3=\{4\}$.
We then use cluster coordinates for $D$, {\em i.e.} apply the linear symplectomorphism 
$\Psi: T^*\bR^{4d}\to  T^*\bR^{4d}$ of the phase space, mapping 
$(p_1,p_2,q_1,q_2)$ to\\
\[(p_D,p^I_D,q_D,q^I_D):={\textstyle \l(p_1+p_2,
\frac {m_2p_1-m_1p_2}{m_D},\frac{m_1q_1+m_2q_2}{m_D},q_1-q_2\ri)}\]
with cluster mass $m_D = m_1+m_2$ and reduced mass
$m^I_D=\frac{m_1 m_2}{m_1+m_2}$ of $D$, and preserving the other coordinates.
In the $C_{1,2}$ center of mass frame, see \eqref{com:c12}
\[\textstyle p_3 = -p_D \qmbox{and} q_3 = -\frac{m_D}{m_3}q_D.\]
We do not rename phase space functions, transformed with $\Psi$.
So in \eqref{Ham:G}
\beq \textstyle
K_{C_1,C_2}= \frac{\|p_D\|^2}{2m_{C_1\!,\!C_2}} \qmbox{,}
H^I_D = \frac{\|p^I_D\|^2}{2m^I_D}+V_{1,2}(q^I_D) 
\Leq{energies}
(whereas $K_D = \frac{\|p_D\|^2}{2m_D}$) and\,%
\footnote{ \ $\frac{m_{C_{1,2}}}{m_3} = \frac{m_D}{m_{C_1,C_2}} > 1$ is the ratio of the distance of the cluster
centers and of $\|q_D\|$.
} 
\beq
\textstyle
V_{C_1,C_2}=
V_{1,3}\big( \frac{m_{C_{1,2}}}{m_3} q_D + \frac{m_2}{m_D} q^I_D \big)+
V_{2,3}\big( \frac{m_{C_{1,2}}}{m_3} q_D -  \frac{m_1}{m_D} q^I_D \big).
\Leq{V:E:C}
Relative angular momentum $L_{C_1,C_2}$ from \eqref{rel:cluster:ang} is a term of
$L_{C_{1,2}}\!=\!\eh\sum_{k=1}^3 q_k\wedge p_k$:
\beq
\textstyle
\hspace*{-2mm}
L_{C_{1,2}}=L_{C_1,C_2}+L^I_D \quad\mbox{with }\ L_{C_1,C_2}= 
\frac{m_D}{m_{C_1\!,\!C_2}}q_D\wedge p_D \ \mbox{ and }\ L^I_D = q^I_D\wedge p^I_D\, .
\Leq{LC12}
Similarly, $J_{C_{1,2}}=\eh\sum_{k=1}^3 m_k\|q_k\|^2$
and the relative moment of inertia \eqref{J:C:D} are related by
\beq
\textstyle
\hspace*{-2mm}
J_{C_{1,2}} =J_{C_1,C_2}+J^I_D\quad\mbox{with } J_{C_1,C_2}=\frac{m_D^2}{2m_{C_1\!,C_2}}\|q_D\|^2
  \mbox{ and } J^I_D=\eh m^I_D \|q^I_D\|^2.
\Leq{JC12}

Assumption \ref{lem:ass:1} in the following lemma is eventually satisfied by {\em all} orbits, for which
asymptotic velocity does not exist. Assumption \ref{lem:ass:2} is satisfied for {\em no} such orbit. 
After showing in Lemma \ref{lem:three:particles:int} 
that Assumption \ref{lem:ass:3} can essentially be skipped, this
will give us an upper bound on the cluster angular momentum of such orbits.
\begin{lemma}[Near-collision of $C_1$ and $C_2$, no interaction with $C_3$] \quad\\
\label{lem:three:particles:no:int}%
Consider for initial condition $\tilde{x}(0)\in \operatorname{NAV}_{E,0}$ the motion 
$t\mapsto \tilde{x}(t)=(\tilde{p}(t),\tilde{q}(t))$ on $\Sigma_{E,0}$,
whose initial condition fulfills the following assumptions:
\begin{enumerate}[1.]
\item 
\label{lem:ass:1}
For some $\delta\in (0,1/2)$, initially the energies are related by (see \eqref{Ham:C:D})
\beq
\textstyle 
|E|\,\le - \delta\, \tilde{H}^I_D(0).
\Leq{tilde:E}
\item 
\label{lem:ass:2}
For the relative cluster angular momentum \eqref{rel:cluster:ang} (with $I$ from Remark \ref{rem:constant:I}), 
\beq
\|\tilde{L}_{C_1,C_2}(0)\| \ge C \ \tilde{K}_{C_1,C_2} (0)^{-\frac{2-\alpha}{2\alpha}} I^{1/\alpha}
( m_{C_1\!,\!C_2} )^{1/2}.
\Leq{l:0}
\item
\label{lem:ass:3}
There is no interaction between the clusters $C_{1,2}$ and $C_3$: 
\[\textstyle
V^E_{\{C_{1,2},C_3\}}(q) \equiv \sum_{i\in C_{1,2},j\in C_3}V_{i,j}(q_i-q_j)=0 \qquad \big(q\in \widehat{M}\big).\]
\end{enumerate}
Then for $C$ in \eqref{l:0} large, one has, for times in $\{t\in \bR \mid \|\tilde{q}_D(t)\| \le 1\}$\,:
\begin{enumerate}[1.]
\item 
\label{tilde:E:L:var}
The variation of relative kinetic energy and cluster angular momentum is small:
\begin{align} 
|\tilde{K}_{C_1,C_2} (t)-\tilde{K}_{C_1,C_2} (0)| &= \cO(C^{-\alpha})\tilde{K}_{C_1,C_2} (0),
\label{rel:kin}\\
\|\tilde{L}_{C_1,C_2}(t)-\tilde{L}_{C_1,C_2}(0)\| &= \cO(C^{-\alpha})\|\tilde{L}_{C_1,C_2}(0)\|,
\label{rel:ang}
\end{align} 
and
\beq 
- (1+2\delta)\tilde{H}^I_D(t)\ \ge\ \tilde{K}_{C_1,C_2}(t).
\Leq{tilde:Hr:interval}
\item 
\label{tilde:qr}
$\|\tilde{q}^I_D \|  < 
\big(\frac{2I}{\tilde{K}_{C_1,C_2} }\big)^{1/\alpha}$\!, but 
$\| \tilde{q}_D\| \ge C\frac{m_{C_1,C_2}}{2m_D}\big(\frac{I}{\tilde{K}_{C_1,C_2} }\big)^{1/\alpha}$, 
thus $\|\tilde{q}^I_D \| \ll \| \tilde{q}_D\|$.
\item 
\label{direction}
The total change of direction $\frac{\tilde{q}'_{C_i}(t)}{\|\tilde{q}'_{C_i}(t)\|}$ 
of the cluster centers is of order 
\beq
 \textstyle\int_\bR \frac{\|\dot{\tilde{q}}_{C_i} \wedge \ddot{\tilde{q}}_{C_i}\|}{\|\dot{\tilde{q}}_{C_i}\|^2}dt
=\cO(C^{-\alpha})\qquad(i=1,2).
\Leq{change:of:dir}
\end{enumerate}
\end{lemma}
\textbf{Proof:}
The proof will be devious, as we first {\em assume} 
statements concerning $\tilde{K}_{C_1,C_2}$, $\tilde{L}_{C_1,C_2}$ and $\tilde{H}^I_D$ that are weaker
than the ones of Assertion \ref{lem:ass:1} (see (\ref{apriori:L}--\ref{apriori:HI})),
then conditionally show all assertions of the lemma, 
and thus prove (in Parts \ref{step3}, \ref{step7} and \ref{step8}) 
that these {\em a priori} inequalites (valid for $t=0$) were justified.

We only consider the three-body case $|C_{1,2}|=3$, since the simpler two-body case $|C_{1,2}|=2$
leads to even better estimates.
Then the nontrivial cluster $D$ equals $C_1$ (or, equivalently, $C_2$).
By Remark \ref{rem:constant:I},
Estimate \eqref{def:moderate} is applicable to the pair potentials $V_{i,j}$ for times $t$ with
$\|\tilde{q}_D(t)\| \le 1$.
\medskip

\noindent
We henceforth assume the {\em a priori} inequalities for all $t$ in $\{t\in \bR \mid \|\tilde{q}_D(t)\| \le 1\}$:
\begin{align}
\label{apriori:K}
|\tilde{K}_{C_1,C_2} (t)-\tilde{K}_{C_1,C_2} (0)| & \le \eh\tilde{K}_{C_1,C_2} (0)\\
\label{apriori:L}
\|\tilde{L}_{C_1,C_2}(t)\| & \ge
\eh C \ \tilde{K}_{C_1,C_2} (t)^{-\frac{2-\alpha}{2\alpha}} I^{1/\alpha}\sqrt{m_{C_1,C_2}}\\ 
\label{apriori:HI}
- 2\tilde{H}^I_D(t)\ &\ge\ \tilde{K}_{C_1,C_2}(t).
\end{align}
The one for $\tilde{L}_{C_1,C_2}(t)$ is similar to \eqref{l:0} for $t=0$. 
For $C$ large, these are weaker than the assertions \eqref{rel:kin}, \eqref{rel:ang} and \eqref{tilde:Hr:interval}.

\begin{enumerate}[1.]
\item
\label{step1}
By \eqref{LC12} in the $C_{1,2}$ center of mass frame 
$L_{C_1,C_2} = \frac{m_D}{m_{C_1\!,\!C_2}} q_D\wedge p_D$,  so that 
by \eqref{energies} and the {\em a priori} inequality \eqref{apriori:L}, 
with $c_1 := \frac{m_{C_1\!,\!C_2}}{2m_D}$,
\beq
\textstyle
\|\tilde{q}_D\| \ge  \frac{m_{C_1\!,\!C_2}}{m_D}  \frac{\|\tilde{L}_{C_1,C_2}\|}{\|\tilde{p}_D\|}
=  \frac{\|\tilde{L}_{C_1,C_2}\|}{m_D} \big(\frac{m_{C_1\!,\!C_2}}{2\tilde{K}_{C_1,C_2}}\big)^{1/2}
\ge  c_1\, C \big(I/\tilde{K}_{C_1,C_2}\big)^{1/\alpha} 
\Leq{q:D:lower}
So conditionally on \eqref{apriori:L}, we have shown the estimate of $\tilde{q}_D$ in Statement~\ref{tilde:qr}.
\item
\label{step2}
The estimate on $\tilde{q}^I_D$ in Statement \ref{tilde:qr}.\ follows from the {\em a priori}
inequality  \eqref{apriori:HI}: 
\beq
\textstyle \|\tilde{q}^I_D\| \le  \big(\frac{I}{|\widetilde{V}_{1,2}|} \big)^{1/\alpha} \le
\big(\frac{I}{|\tilde{H}^I_D|} \big)^{1/\alpha}  \le
\big(\frac{2I}{\tilde{K}_{C_1,C_2}} \big)^{1/\alpha}\, .
\Leq{q:ID:upper}

So \eqref{apriori:HI} implies for $C$ large that the relative cluster energy
$\tilde{H}_{C_1,C_2} (t)$ is mainly kinetic in the following sense:
\beq
0<I\|\tilde{q}_D\|^{-\alpha}\leq c_1^{-\alpha} C^{-\alpha} \tilde{K}_{C_1,C_2}\,  .
\Leq{mainly:kinetic}
In particular, for $c_2:= 3c_1^{-\alpha} $ and $C$ large,
\beq\textstyle
\hspace*{-3mm}
|\tilde{V}_{C_1,C_2}|\le c_2\, C^{-\alpha} \tilde{K}_{C_1,C_2}\qmbox{and }
|\langle q_D, \nabla_{q_D} \tilde{V}_{C_1,C_2}\rangle |\le c_2\, C^{-\alpha} \tilde{K}_{C_1,C_2}.
\Leq{ll:ll}
\item
\label{step3}
As by Assumptions \ref{lem:ass:1}.\ and \ref{lem:ass:3}., the constant $\tilde{H}^I_{C_{1,2}}$ is bounded
from above by 
\[\tilde{H}^I_{C_{1,2}}
=E-\tilde{H}_{\{C_{1,2},C_3\}}=E-\tilde{K}_{\{C_{1,2},C_3\}}\le E\le -\delta\, \tilde{H}^I_D(0),\] 
and by \eqref{ll:ll} and \eqref{Ham:G}, 
\[0<(1-c_2\, C^{-\alpha})\tilde{K}_{C_1,C_2}\le \tilde{H}_{C_1,C_2}
= \tilde{H}^I_{C_{1,2}} - \tilde{H}^I_{D} \le -(1+\delta) \tilde{H}^I_{D}.\]
So $ -(1+2\delta) \tilde{H}^I_{D}\ge \tilde{K}_{C_1,C_2}$ for $C$ large, showing \eqref{tilde:Hr:interval}
conditionally on (\ref{apriori:L}--\ref{apriori:HI}).
\item
\label{step4}
To control the time evolution of $q_D$, we consider $\tilde{J}_{C_1,C_2}$ 
in the $C_{1,2}$--center of mass frame, see \eqref{JC12}. 
Then 
\beq
\textstyle
\frac{d^2}{dt^2} \tilde{J}_{C_1,C_2} = 
\textstyle 2\tilde{K}_{C_1,C_2}+\langle \tilde{q}_D, \nabla_{q_D} V_{C_1,C_2}(\tilde{q})\rangle 
\Leq{J:C:pp}
Then with \eqref{ll:ll}, for $C$ large \eqref{J:C:pp} leads to the inequality
\[\textstyle
\frac{d^2}{dt^2} \tilde{J}_{C_1,C_2} \ge \tilde{K}_{C_1,C_2} .\]
\item
\label{step5}
So with the {\em a priori} estimate \eqref{apriori:K} 
(implying $\tilde{K}_{C_1,C_2} (t)\ge \eh \tilde{K}_{C_1,C_2} (0)>0$),
\[\textstyle
\frac{d^2}{dt^2} \tilde{J}_{C_1,C_2} \ge \eh \tilde{K}_{C_1,C_2} (0)>0.\] 
Thus there is a unique time $t_0$ where $ \tilde{J}_{C_1,C_2}$ attains its minimum, and 
by \eqref{JC12}, \eqref{q:D:lower} and \eqref{apriori:K} we get the propagation estimate
\begin{align}
\tilde{J}_{C_1,C_2}(t) 
& \ge \tilde{J}_{C_1,C_2}(t_0) + {\textstyle \frac{1}{4}}\tilde{K}_{C_1,C_2} (0)\, (t-t_0)^2 \NN\\
& \ge \textstyle  \frac{1}{8}m_{C_1\!,C_2}\, C^2 \,(\frac{2}{3}I/\tilde{K}_{C_1,C_2} (0)))^{\frac{2}{\alpha}} 
+ \frac{1}{4} \tilde{K}_{C_1,C_2} (0)\, (t-t_0)^2.
\label{prop:est}
\end{align}
So in both time directions the
distance between the cluster centers $q_{C_1}$ and $q_{C_2}$ diverges (at least) linearly.
\item
\label{step6}
We now show Statement \ref{direction}., estimating $\int_\bR  (\tilde{J}_{C_1,C_2})^{-(1+\alpha)/2} \,dt$ 
by use of \eqref{prop:est} and with
\beq
\textstyle
\int_\bR \left(b+c t^2\right)^{-(1+\alpha)/2}\,dt = 
\frac{C_1(\alpha)} {b^{\alpha/2} \sqrt{c} }\qmbox{for}
C_1(\alpha):=\frac{\sqrt{\pi }\Gamma \left(\frac{\alpha}{2}\right)}{\Gamma \left(\frac{\alpha+1}{2}\right)}.
\Leq{integration:const}
Thus
\beq
\int_\bR  (\tilde{J}_{C_1,C_2})^{-(1+\alpha)/2} \,dt\le \textstyle
24\, C_1(\alpha) m_{C_1,C_2}^{-\alpha/2}\; C^{-\alpha}\frac{\sqrt{\tilde{K}_{C_1,C_2} (0)}}{I}.
\Leq{J:integral}
Using \eqref{JC12} and \eqref{J:integral}, the total change of direction is bounded from above by 
\begin{align}
\textstyle \int_\bR \frac{\|\dot{\tilde{q}}_D \wedge \ddot{\tilde{q}}_D\|}{\|\dot{\tilde{q}}_D\|^2}\, dt
&\textstyle\le \int_\bR \frac{\|\ddot{\tilde{q}}_D\|}{\|\dot{\tilde{q}}_D\|}\, dt
\le \frac{\sqrt{m_{C_1,C_2}}}{m_D}\int_\bR \frac{\|\nabla \tilde{V}_{C_1,C_2}\|}{\sqrt{\tilde{K}_{C_1,C_2} }} \, dt
\label{change:of:direction}\\
&\textstyle\le \frac{\sqrt{2m_{C_1,C_2}/\tilde{K}_{C_1,C_2} (0)}}{m_D}\int_\bR \|\nabla \tilde{V}_{C_1,C_2}\|\,dt
\NN\\
&\textstyle
\le (1+c)\alpha  \frac{I}{m_D} \sqrt{\frac{2m_{C_1,C_2}}{\tilde{K}_{C_1,C_2} (0)}}  
\int_\bR \|\tilde{q}_D\|^{-1-\alpha}\,dt 
\NN\\
&\textstyle\le 2 \frac{I}{m_D} \sqrt{\frac{2m_{C_1,C_2}}{\tilde{K}_{C_1,C_2} (0)}} 
(\frac{m_D}{\sqrt{2m_{C_1,C_2}}})^{\alpha+1}\int_\bR  (\tilde{J}_{C_1,C_2})^{-(1+\alpha)/2} \,dt \NN\\
&\textstyle\le 48\,C_1(\alpha)\,
(\frac{m_D}{m_{C_1,C_2}})^{\alpha}\ C^{-\alpha}\, \NN
\end{align} 
for $C$ large. So as $C\nearrow \infty$, the total change of direction of $\tilde{q}'_{C_1}=\tilde{q}'_{D}$ 
goes to zero.
\eqref{change:of:dir} follows for cluster $C_2$, too, since 
$\frac{\tilde{q}'_{C_2}(t)}{\|\tilde{q}'_{C_2}(t)\|}=-\frac{\tilde{q}'_{C_1}(t)}{\|\tilde{q}'_{C_1}(t)\|}$
in the center of mass frame.
\item
\label{step7} 
We next estimate the time evolutions of the energy terms, using the propagation estimate
\eqref{prop:est}. We begin with $\tilde{H}^I_D$, see \eqref{energies}.
A priori, the time derivative
\[\textstyle
\frac{d}{dt} \tilde{H}^I_D = \langle \tilde{v}^I_D, \nabla_{\! q^I_D} V_{C_1,C_2}(\tilde{q})\rangle\]
of the internal energy of cluster $D$ could be unbounded, since near-collisions of the particles 1 and 2
can lead to large relative velocities $\tilde{v}^I_D$. Therefore we apply partial integration when estimating
the change of $\tilde{H}^I_D$.
\begin{align*}
\textstyle
\tilde{H}^I_D(t)-\tilde{H}^I_D(0) &\textstyle
= \int_0^t \frac{d}{dt} \tilde{H}^I_D(s)\,ds\\
&\textstyle
= \langle \tilde{q}^I_D, \nabla_{\! q^I_D} V_{C_1,C_2}(\tilde{q})\rangle\big|_0^t - 
\int_0^t \langle \tilde{q}^I_D, \frac{d}{ds} \nabla_{\! q^I_D} V_{C_1,C_2}(\tilde{q})\rangle\,ds.
\end{align*}
\begin{enumerate}[$\bullet$]
\item 
By Statement \ref{tilde:qr}, $\|\tilde{q}^I_D\| \ll \|\tilde{q}_D\|$ if $C$ is large.
Then using \eqref{mainly:kinetic} and the a priori inequality,
the first term is bounded from above by $6c_1^{-\alpha} C^{-\alpha} \tilde{K}_{C_1,C_2}(0)$. 
\item 
Similarly we bound the integrand of the second term by
\[\textstyle
|\langle \tilde{q}^I_D, \frac{d}{ds} \nabla_{\! q^I_D} V_{C_1,C_2}(\tilde{q})(s)|
\le 4I\|\tilde{q}(s)\|^{-1-\alpha}\sqrt{\tilde{K}_{C_1,C_2}(0)/m_{C_1,C_2}}.\]
Integration using \eqref{J:integral} leads to the bound
\[\textstyle
|\int_0^t \langle \tilde{q}^I_D, \frac{d}{ds} \nabla_{\! q^I_D} V_{C_1,C_2}(\tilde{q})\rangle\,ds|
\ \le\ 68\, C_1(\alpha) \; C^{-\alpha} \tilde{K}_{C_1,C_2}(0).
\]
\end{enumerate}
So the variation of the internal energy of cluster $D$ is bounded by
\beq
\textstyle
|\tilde{H}^I_D(t)-\tilde{H}^I_D(0)|\le\ 74\, C_1(\alpha) \; C^{-\alpha} \tilde{K}_{C_1,C_2}(0).
\Leq{H:var}
For $C$ large this is stronger than the {\em a priori\,} inequality \eqref{apriori:HI}.\\
The same bound applies to the variation of $\tilde{H}_{C_1,C_2}$, as the internal energy 
$H^I_{C_{1,2}}=H_{C_1,C_2}+H^I_{D}$ of $C_{1,2}$ is constant, see \eqref{Ham:G}.

The estimate for the time evolution of $\tilde{V}_{C_1,C_2}$ follows directly from \eqref{ll:ll}:
\beq
|\tilde{V}_{C_1,C_2}(t) - \tilde{V}_{C_1,C_2}(0)|\le 2c_2\, C^{-\alpha} \tilde{K}_{C_1,C_2}.
\Leq{V:var}
Def.\ \eqref{Ham:C:D} finally shows, using \eqref{H:var} and \eqref{V:var}, 
that with $c_3:=2c_2+74C_1(\alpha)$
\beq
|\tilde{K}_{C_1,C_2}(t) - \tilde{K}_{C_1,C_2}(0)|\le c_3 \, C^{-\alpha} \tilde{K}_{C_1,C_2}.
\Leq{K:var}
For $C$ large this is stronger than the {\em a priori} inequality \eqref{apriori:K}
and will show assertion \eqref{rel:kin}.
\item
\label{step8}
The time evolution of the cluster angular momentum follows from
\[\textstyle 
\frac{d}{dt}\tilde{L}_{C_1,C_2}=-\frac{m_D}{m_{C_1,C_2}}
\tilde{q}_D\wedge \nabla _{q^I_D} V_{C_1,C_2}(\tilde{q}).\]
So by an estimate similar to \eqref{ll:ll} 
\begin{align}
\label{ang:mom:ansatz}
\textstyle \Big|\frac{d}{dt}\|\tilde{L}_{C_1,C_2}\| \Big| &\le \textstyle 
\big\|\frac{d}{dt}\tilde{L}_{C_1,C_2}\big\|
\le 2\alpha m_Dm_3\frac{\|\tilde{q}^I_D\|}{\big(\frac{m_D}{m_{C_1,C_2}}\|\tilde{q}_D\|\big)^{1+\alpha}} \\
&\textstyle 
=  \alpha m_Dm_3 (\frac{m_{C_1,C_2}}{2})^{(1+\alpha)/2} \|\tilde{q}^I_D\| (\tilde{J}_{C_1,C_2})^{-(1+\alpha)/2} \NN
\end{align} 
By \eqref{J:integral} 
and Statement \ref{tilde:qr}.\ on $\|\tilde{q}^I_D\|$
\begin{align}
\|\tilde{L}_{C_1,C_2}(t)-\tilde{L}_{C_1,C_2}(0)\| 
&\textstyle\le 24\,C_1(\alpha)
\sqrt{m_{C_1,C_2}} \|\tilde{q}^I_D\| C^{-\alpha}\,\sqrt{\tilde{K}_{C_1,C_2} (0)}, \NN\\
&\textstyle\le C_3 \ 
C^{-\alpha}\ \tilde{K}_{C_1,C_2} (0)^{-\frac{2-\alpha}{2\alpha}}I^{1/\alpha} \sqrt{m_{C_1,C_2}} , \NN
\end{align} 
with $C_3:= 24\,C_1(\alpha) \big(\frac{m^I_D}{m_3} \big)^{1/\alpha}$.
Comparing with \eqref{l:0}, the relative variation~of cluster angular momentum goes to zero
as $C\nearrow \infty$, in accordance with {\em a priori} inequalities (\ref{apriori:L}-\ref{apriori:HI})
and assertion \eqref{rel:ang}. This finishes the proof.
\hfill~$\Box$
\end{enumerate}
We now consider the influence of cluster $C_3$.
\begin{lemma}[Near-collision of $C_1$ and $C_2$, interacting with $C_3$] \quad\\
\label{lem:three:particles:int}%
If $\|\tilde{q}_{C_3}(t)-\tilde{q}_{C_1}(t)\| \gg 1$ for $\{t\in \bR \mid \|\tilde{q}_D(t)\| \le 1\}$, 
the conclusions of Lemma \ref{lem:three:particles:no:int} are still valid (with slightly worse constants),
without its Assumption~\ref{lem:ass:3}.
\end{lemma}
\textbf{Proof:}
This follows from a perturbation argument concerning the forces, as in the proof of Lemma
\ref{lem:three:particles:no:int} we used the estimate \eqref{def:moderate}
($|\pa_\beta V_{i,j}(q)| \le  I \|q\|^{-\alpha - |\beta|}$ for $\|q\|$ small, see also Remark \ref{rem:constant:I}).
If $\|\tilde{q}_{C_3}(t)-\tilde{q}_{C_1}(t)\| \gg 1$ but $\|\tilde{q}_D(t)\| \le 1$, then 
\hfill $\Box$\\[2mm]
So in particular, for $x\in\operatorname{NAV}_{E,0}$ 
one has the reverse estimate for the relative angular momentum \ref{rel:cluster:ang}
\beq
\|\tilde{L}_{C_1,C_2}(t)\| \le 
2 C \ \tilde{H}_{C_1,C_2} (0)^{-\frac{2-\alpha}{2\alpha}} I^{1/\alpha}\sqrt{m_{C_1,C_2}}
\qquad(t\in [0,t_{\max}]).
\Leq{reverse:NAV}
This will be used when we show that the orbit hits almost all Poincar\'e surfaces.
\subsection{Motion of the messenger cluster} \label{sub:sec:prop:est}
%
We consider the following setting: For $(C_1,C_2,C_3)\in \cT$ the cluster center $q_{C_2}$ of $C_2$
moves in a time interval $(s_1,s_3)$ from a neighborhood of $q_{C_1}$ to the one of $q_{C_3}$. 
The respective neighborhoods are defined, using energy considerations.\\
For the total energy $E$ we can assume, using Lemma \ref{lem:nontrivial:cluster} (c), 
that the internal energy $\tilde{H}^I_D(s_1)$ of the 
nontrivial cluster $D$ is so negative that
$|E| \ll -\tilde{H}^I_D(s_1)$. The times $s_1$ and $s_3$ are chosen so that
the initial distances between $q_{C_2}$ and $q_{C_i}$
satisfy the inequality 
\beq
\| \tilde{q}_{C_2}-\tilde{q}_{C_i} \|(s_i)\ge 
C\, \max\big( (\tilde{K}_{C_2,C_i})^{-1/\vep}(s_i)\, , \, 1\big) \qquad (i=1,3).
\Leq{initial}
This condition means that in the time interval $(s_1,s_3)$, for $C$ large enough
\begin{enumerate}[1.]
\item 
one can apply the long range estimate \eqref{long:range}, in particular 
$V_{i,j}(q)=\cO(\|q\|^{-\vep})$, if $i$ and $j$ belong to different clusters;
\item 
thus the distances of the clusters are so large that the energies
$\tilde{H}_{C_2,C_i}$ are mainly kinetic: $|\tilde{V}_{C_2,C_i}(t)|\ll \tilde{K}_{C_2,C_i}(t)$ for 
$t\in (s_1,s_3)$ and $C\gg1$ in \eqref{initial}.
\end{enumerate}

\begin{lemma}[Variation of angular momenta and directions]
\label{lem:small:var:a:m}%
With \eqref{initial}, 
\begin{enumerate}[1.]
\item 
the variations of the relative angular momenta are of order
\beq
\|\tilde{L}_{C_2,C_i}(s_1)-\tilde{L}_{C_2,C_i}(s_3)\|
= \cO \big( \big(\tilde{K}_{C_1,C_2} \big)^{-1/\alpha}(s_1) \big) \qquad (i=1,3),
\Leq{L:L:s13}
\item 
the variations of the directions of the three clusters are of order
\beq
\scalebox{1.3}{$\int$}_{\!\!s_1}^{s_3}\, \textstyle  \frac{\|\tilde{q}'_{C_2,C_i} \wedge\,\tilde{q}''_{C_2,C_i}\|}
{\|\tilde{q}'_{C_2,C_i}\|^2}\, dt
= \cO \big( \big(\tilde{K}_{C_1,C_2} \big)^{-1/\alpha}(s_1) \big) \qquad (i=1,3),
\Leq{cod:messenger}
\item 
Consider the affine line $\{ \tilde{q}_{C_2}(s_2)\} + \operatorname{span}(\tilde{p}_{C_2}(s_2)) $, 
with $s_2:=\eh(s_1+s_3)$. Then the maximal distance of $\tilde{q}_{C_2}(t)$ for $t\in (s_1,s_3)$
is of order $\cO \big( \big(\tilde{K}_{C_1,C_2} \big)^{-1/\alpha}(s_1) \big)$. Similar statements hold for
$C_1$ and $C_3$.
\end{enumerate}
Thus all vanish in the limit $t\nearrow T(x)$.
\end{lemma}
\begin{remark} 
Note that for $i=3$, too the order depends on $\tilde{K}_{C_1,C_2}(s_1)$.
This is because the kinetic energy comes into play through the estimate 
$\|\tilde{q}^I_D \|  < \big(\frac{2I}{\tilde{K}_{C_1,C_2} }\big)^{1/\alpha}$ in 
Lemma \ref{lem:three:particles:no:int}  for the size of the nontrivial cluster $D$.
The difference is important, since 
unlike $\tilde{K}_{C_2,C_3}(s_3)$, $\tilde{K}_{C_1,C_2}(s_1)$ is bounded below,
by $\tilde{K}^E_{C_2}$ from Lemma \ref{lem:part:energy}, since after collision $C_1$ and $C_3$ 
move in nearly opposite directions.~$\Diamond$
\end{remark}
\textbf{Proof:}
We proceed like in the proof of Lemma \ref{lem:three:particles:no:int}, with less details.\\[-6mm]
\begin{enumerate}[1.]
\item 
We use the inequalities
$\tilde{J}_{C_2,C_i}(t) \ge\tilde{J}_{C_2,C_i}(s_i) +  \frac{1}{4} \tilde{K}_{C_2,C_i} (s_i)\, (t-s_i)^2$,
for $t\in(s_1,s_2)$, similar to \eqref{prop:est}. So 
\begin{align}\textstyle
\int_{s_1}^{s_3}   
(\tilde{J}_{C_2,C_i})^{-(1+\vep)/2} \,dt&\le\textstyle
\frac{2\, 2^{(1+\vep)/2} }{(\tilde{K}_{C_2,C_i} (s_i))^{1/2}}
\int_{(\tilde{J}_{C_2,C_i} (s_i))^{1/2}}^\infty t^{-1-\vep}\, dt
\NN \\
&=\textstyle
\frac{2\, 2^{(1+\vep)/2} }{\vep (\tilde{K}_{C_2,C_i} (s_i))^{1/2}}(\tilde{J}_{C_2,C_i} (s_i))^{-\vep/2} \, ,
\label{infinite:J:int}
\end{align}
applying Cauchy-Schwarz. 
Thus by \eqref{initial}, the is integral bounded independent of $\tilde{K}_{C_2,C_i} (s_i)>0$, and of
order $\cO\big( (\tilde{K}_{C_2,C_i} (s_i))^{-1/2}\big)$ for large relative kinetic cluster energies.
With $ \big|\frac{d}{dt}\|\tilde{L}_{C_2,C_i}\| \big| 
\le  c \|\tilde{q}^I_D\| (\tilde{J}_{C_2,C_i})^{-(1+\alpha)/2} $, see
\eqref{ang:mom:ansatz}, and with \eqref{infinite:J:int},  \eqref{L:L:s13} follows, more precisely:
\[\textstyle
\|\tilde{L}_{C_2,C_i}(s_1)-\tilde{L}_{C_2,C_i}(s_3)\|
= \cO\big(\min\big(\tilde{K}_{C_2,C_i} (s_i))^{-1/2},\,1\big)\ \tilde{K}_{C_1,C_2} (s_1)^{-\frac{1}{\alpha}}\big).
\]
\item 
Unlike in \eqref{change:of:direction}, we use that the relative velocities of the clusters are nearly parallel
to their relative positions (and thus their accelerations, as the $V_{i,j}$ are central). 
Thus we use a first order Taylor estimate for the integrand in \eqref{cod:messenger}.
By \eqref{l:0} (remember that the $\operatorname{NAV}$ satisfy the {\em converse} inequality), initially 
\[ \textstyle
\frac{\|\tilde{p}_{C_2,C_i} (t) \wedge \tilde{q}_{C_2,C_i} (t)\|}{\| \tilde{p}_{C_2,C_i} (t)\|} 
= \frac{\| \tilde{L}_{C_2,C_i}(t)\|}{\|\tilde{p}_{C_2,C_i}(t)\|}
= \cO \big(\tilde{K}_{C_1,C_2} (t)^{-\frac{1}{\alpha}} \big)\qquad \big(t\in(s_1,s_3)\big),\]
which is the same order as the size $\|\tilde{q}^I_D(t)\|$ of the nontrivial cluster. Therefore,
\begin{align} \textstyle
\frac{\|\tilde{q}'_{C_2,C_i} \wedge\,\tilde{q}''_{C_2,C_i}\|}{\|\tilde{q}'_{C_2,C_i}\|^2}(t)
&\textstyle 
=  \cO \Big(\frac{\|D^2\tilde{V}_{C_2,C_i}(t)\|}{\|\tilde{q}'_{C_2,C_i}(t)\| \|\tilde{q}_{C_2,C_i} (t)\|}
 \ \tilde{K}_{C_1,C_2} (t)^{-\frac{1}{\alpha}} \Big)
 \label{curvature}\\
&\textstyle
= \cO \big((\tilde{J}_{C_2,C_i})^{-(3+\vep)/2}(t) \ \tilde{K}_{C_1,C_2} (t)^{-\frac{1}{\alpha}} \big)
\qquad \big(t\in(s_1,s_3)\big)\NN
\end{align}
Integration leads to the estimate
\begin{align*} \textstyle
\scalebox{1.3}{$\int$}_{\!\!s_1}^{s_3}\, \textstyle  \frac{\|\tilde{q}'_{C_2,C_i} \wedge\,\tilde{q}''_{C_2,C_i}\|}
{\|\tilde{q}'_{C_2,C_i}\|^2}\, dt
&= \textstyle\cO\Big(
\int_{s_i}^\infty   t^{-(3+\vep)/2} \,dt \,K_{C_2,C_i}(s_i)^{-(4+\vep)/2}\tilde{K}_{C_1,C_2} (s_i)^{-\frac{1}{\alpha}} 
\Big)\\
&=
\cO \big(\tilde{K}_{C_1,C_2} (s_1)^{-\frac{1}{\alpha}} \big). 
\end{align*}
\item
The third statement follows by double integration of the curvature \ref{curvature} for the space curve 
$\tilde{q}_{C_2}$, starting at time $s_2$.\hfill $\Box$
\end{enumerate}
%
\section{Applicability of the Poincar\'e surface method} \label{sect5}
%
To prove the main theorem (page \pageref{thm:main}), we show  in Section \ref{sub:sect:NAV:wanders}
that $\operatorname{NAV}\subseteq\operatorname{Wand}$,
and in Section \ref{sub:sect:NAV:transitional} that 
$\operatorname{NAV}\subseteq \operatorname{Trans}$.
Then the Poincar\'e surface method of Theorem \ref{thm:FK} is applicable.
%
\subsection{\protect{$\operatorname{NAV}$} is wandering} \label{sub:sect:NAV:wanders}
%
\textbf{Proof of Theorem \ref{thm:NAV:wandering}:}
The proof is based on the extension of von Zeipel's theorem in \cite[Chapter 2.4]{Fl}, see also
\cite[Chapter 12.6]{Kn}. 
As for lack of space we cannot fully reproduce that proof here, we indicate the main points.

It uses the Graf partition \eqref{Z4} of configuration space, originally devised in the context of quantum 
$n$-body scattering. This partition and its associated convex function \eqref{J:delta} allow to 
focus attention on the motion of the cluster centers, instead of the more complicated 
cluster-internal motion, see 
\cite[Chapter 5]{DG}.

\begin{figure}[htbp]
\vspace*{-3mm}
\begin{center}
\includegraphics[height=50mm]{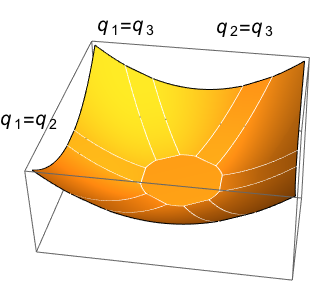}
\caption{The function $J^{(\delta)}:M\to \bR$ on the $n$-particle configuration space $M$, 
for $n=3$ particles in $d=1$ dimension and  in the center of mass frame}
\label{fig:J:delta}
\end{center}
\end{figure}

For the partition lattice $\cP(N)$ of $N=\{1,\ldots,n\}$, $\cC\in\cP(N)$,
the projections $\Pi_{\cC}^E, \Pi_{\cC}^I:M\to M$ to the cluster-external/internal coordinates
(see \eqref{eq:DefDeltaE}), one sets
\beq
J=J_\cC^E+J_\cC^I \qmbox{with} 
J_\cC^E := J\circ\Pi_\cC^E \qmbox{and} J_\cC^I := J\circ\Pi_\cC^I\,,
\Leq{Z8}
and, see Figure \ref{fig:J:delta}, for $\delta\in (0,1]$ sufficiently small
\beq
J^{(\delta)}:M\to\bR\qmbox{,} J^{(\delta)}(q)
:= \max \big\{ J_\cC^E(q) + \delta^{|\cC|} \;\big|\;\cC\in\cP(N) \big\}.
\Leq{J:delta}

\noindent
The {\em Graf partition} of the configuration space  $M$ is the family of subsets 
\beq
\Xi_\cC^{(\delta)}:=
\l\{q\in M\;\big|\;J_{C_1,C_2}(q)+\delta^{|\cC|} = J^{(\delta)}(q) \ri\}
\qquad \bigl(\cC\in\cP(N)\bigr).
\Leq{Z4}
$\Xi_\cC^{(\delta)}$ should be compared with $\Xi_\cC^{(0)} $ from \eqref{Xi:C:0}.\\ 
Omitting the initial condition $x\in P$, there are piecewise constant mappings 
\beq \textstyle
\cA:\big(0,T(x)\big) \to\cP(N) \qmbox{with}
\frac{\tilde{q}(t)}{t^{1-\vep/2}}\in\Xi_{\cA(t)}^{(\delta)}\, 
\Leq{00V}
(unique up to points on boundaries between the $\Xi_{\cA(t)}^{(\delta)}$)
and thus  for $\tilde{Q}(t):=\frac{\tilde{q}(t)}{t}$ 
the continuous, piecewise differentiable approximant of $j$ (see \eqref{j:x}), 
\[\tilde{j}^{(\delta)}(t) := \tilde{j}^E(t)+\delta^{|\cA(t)|}t^{-\vep}
\qmbox{with}
\tilde{j}^E(t) := J_{\cA(t)}^E\big( \tilde{Q}(t)\big) 
= \eh \langle \tilde{Q}^E_{\cA(t)} , \tilde{Q}^E_{\cA(t)} \rangle_\cM.\] 
Omitting $\cA$, the time derivative  
\beq
\textstyle \frac{d}{dt}\tilde{j}^E = \frac{1}{2t} \langle \tilde{Q}^E ,\frac{d}{dt}\tilde{q}^E - \tilde{Q}^E \rangle_\cM 
\Leq{djE}
of the external part $\tilde{j}^E$ is written as a sum of a nonnegative function
and a function on $(t_0,T(x))$ whose modulus is smaller than $C\,t^{-1-\epsilon'}$, 
for some $\epsilon'\in(0,\vep)$ and for some $C>0$, only depending on $V$ and not on $x$. 
This shows existence of 
\[ \textstyle \lim_{t\nearrow T(x)} \tilde{j}^E(t) \,\in\, [0,\infty].\]
\begin{enumerate}[1.]
\item 
The remarks after \eqref{djE} show that 
\beq
\textstyle  
\|\frac{d}{dt}\tilde{q}^E(t) \|_\cM \ge 
\frac{\langle \tilde{Q}^E(t) ,\frac{d}{dt}\tilde{q}^E(t) \rangle_\cM}{\|\tilde{Q}^E(t)\|_\cM}
\ge \|\tilde{Q}^E(t)\|_\cM-\cO(t^{-\epsilon'}). 
\Leq{speed}
For $x\in \operatorname{NAV}$ the term $ \|\tilde{Q}^E(t)\|_\cM=(2\tilde{j}^E(t))^{1/2}$ on the right side
diverges: $\lim_{t\nearrow T(x)}\|\tilde{Q}^E(t)\|_\cM=\infty$. 
Thus the external speed $\|\frac{d}{dt}\tilde{q}^E(t) \|_\cM$ has the same limit. 
So kinetic energy goes to $\infty$, too and $\lim_{t\nearrow T(x)} \tilde{V}(t)=-\infty$.
\item 
As escape time $T:P\to (0,\infty]$ is lower semicontinuous, 
for any $\tau\in (0, T(x))$ there is a neighborhood $U$
of $x$ so that $T|_U\ge \tau$. By shrinking $U$, if necessary, from continuity of the flow $\Phi$
and \eqref{speed} we get for $y\in U$
\[ j^E\!\circ\Phi_\tau(y)\ge j^E\!\circ\Phi_\tau(x)-1\ \mbox{ and }\ 
j^E\!\circ\Phi_t(y) \ge  j^E\!\circ\Phi_\tau(y)-C t^{-\epsilon'} 
\ \ \big(t\in \big[\tau, T(y)\big)\big).  \]
As by increasing $\tau$, $j^E\!\circ\Phi_\tau(x)$ can be chosen arbitrarily large,
$j^+$ is lower semicontinuous at $x\in \operatorname{NAV}$.
But as $j^+(x)=\infty$, $j^+$ is continuous at~$x$.
\item 
As $J^E$ is bounded on $U$, this also shows that $x_0$ is wandering.
As $j^+$ is lower semicontinuous at $\operatorname{NAV}$, 
for $k\in \bN$ there are open neighborhoods $\operatorname{NAV}_k\supseteq \operatorname{NAV}$
with $j^+(y)>k$ for $y\in\operatorname{NAV}_k$. 
So $\operatorname{NAV}=\bigcap_{k\in\bN}\operatorname{NAV}_k$ is Borel.
\hfill $\Box$
\end{enumerate}
%

\subsection{\protect{$\operatorname{NAV}$} is transitional} \label{sub:sect:NAV:transitional} 

For the admissible potentials $V$ of Definition \ref{admissible}, the pair potentials are central. 
So the total angular momentum $L$, see \eqref{total:ang:momentum}, is a constant of the motion.

This is also the case in the nondeterministic model of Sect.\ \ref{sect3}. But there, additionally the 
values of the cluster angular momenta are bounded by $|L|$, see~\eqref{Li:bound}. 

Although that last property is not expected for the true dynamics and all times $t$, 
that is the case in the limit $t\nearrow T(x)$:
\begin{lemma} [Boundedness of the cluster angular momenta] \quad\\ 
\label{lem:cluster:ang:mom}
For all $\ell>0$ and $x\in \operatorname{NAV}$ there is a minimal time $t_0\in (0,T(x))$ with
\[ \| \tilde{L}^E_C(t) \| \le \|L(x)\| + \ell \qmbox{and} \|\tilde{L}^I_C(t)\| \le \ell 
\qquad \big(t\in [t_0, T(x)) , C\in\cA(t)\big). \] 
Concerning the $\tilde{L}^E_C$ estimate, the term $\ell>0$ is only needed if $\,\|L(x)\|=0$.
\end{lemma}
\textbf{Proof:}
We have, see \eqref{ex:ang:mom} and \eqref{rel:ang:mom}, the relation
\[\textstyle
\tilde{L}(t) = \sum_{C\in \cA(t)} \big(\tilde{L}^E_C(t) +\tilde{L}^I_C(t) \big) \qquad \big(t\in \big[0,T(x)\big)\big), \]
and know from Section \ref{sub:sect:B2} that for $t\in [t_0, T(x))$, the number of clusters is
$|\cA(t)| \in \{2,3\}$, changing infinitely often between the two cases. 
\begin{enumerate}[1.]
\item 
For $t\in  (t_{2\ell},t_{2\ell+1})$, $|\cA(t)| = 3$ and there is a unique cluster
$D(t)\in \cA(t)$ of size $|D|=2$. 
By Lemma \ref{lem:nontrivial:cluster} (c) its internal energy $\tilde{H}^I_{D}=\tilde{K}^I_{D}+\tilde{V}^I_{D}$ 
goes to $-\infty$, and so does $\tilde{V}^I_{D}\le \tilde{H}^I_{D}$.
Using admissibility of $V$, 
$\|\tilde{q}^I_{D}\| = \cO(|\tilde{V}^I_{D}|^{-1/\alpha})= \cO(|\tilde{H}^I_{D}|^{-1/\alpha})$, whereas
$\|\tilde{p}^I_{D}\| = \cO(|\tilde{H}^I_{D}|^{-1/2})$.
So internal angular momentum 
\[ \|\tilde{L}^I_{D(t)}(t) \| = \cO\big((-\tilde{H}^I_{D(t)} (t))^{^{-\frac{2-\alpha}{2\alpha}}} \big)
\qquad \big(t\in\cup_{\ell\in\bN} (t_{2\ell},t_{2\ell+1})\big),\] 
going to zero as $t\nearrow T(x)$.
$t$~belongs to an interval $(t_{2\ell},t_{2\ell+1})$ on which $\cA$ is constant.
Lemma \ref{lem:small:var:a:m} shows that on this interval the variation of the $\tilde{L}_C$ goes to zero, 
and at the end points, the relative angular momenta of the colliding clusters
converge to zero as $t\nearrow T(x)$. 
In the kinematical formulae \eqref{3L:L}, $|L_C| \le (1-\frac {m_{\min}}{M}) |L|$. 
As $\frac {m_{\min}}{M}\in (0,1)$,  adding these dynamical terms
finishes the proof for $|\cA(t)| = 3$
by increasing $t_0$ appropriately.
\item 
For $t\in (t_k,t_{k+1})$ with $k$ odd, that is $|\cA(t_k^+)|=|\cA(t_{k+1}^-)| = 2$ clusters, 
one already knows from part 1.\ of the proof that the three clusters present for $t\nearrow t_k$
respectively $t\searrow t_{k+1}$ meet the bounds of the lemma. 
A propagation estimate like \eqref{prop:est} shows that within the interval $ (t_k,t_{k+1})$ the weak interaction
between the two clusters does not suffice to violate these bounds.~$\quad\Box$
\end{enumerate}
So we can exhaust the set $\operatorname{NAV}_{E,0}$ by first considering only the subset
whose cluster angular momenta are bounded by  ${\frak L}>0$, show that this subset is of Liouville
zero, and then let  ${\frak L}\nearrow \infty$.
 
Usually omitting the indices $E$ and ${\frak L}$, for a given value ${\frak L}>0$ 
we now define the family of Poincar\'e surfaces 
\beq
\cH_m \equiv \cH_{m,E,{\frak L}}  
:= \bigcup_{\bar{\cC}\in {\cal Q}} \cH_m^{\bar{\cC}} \subseteq \Sigma_{E,0}\qquad(m\in \bN),
\Leq{the:surfaces}
with $\bar{\cC}=(C_1,C_2,C_3)$ and $\{C_1,C_2,C_3\}\in\cP_3(N)$, see page \pageref{combinatorics}.
$\cH_m^{\bar{\cC}}$ is defined so that the center of mass $q_{C_2}$ 
of the messenger cluster $C_2$ is in the hyperplane 
\beq
S_m(q_{C_1}) := \textstyle
\big\{r\in \bR^d\mid \big\langle r,\frac{q_{C_1}}{\|q_{C_1}\|}\big\rangle = \|q_{C_1}\|-1/m \big\} 
\Leq{def:hyper}
perpendicular to $q_{C_1}$,
that has minimal distance $1/m$ from $q_{C_1}$.
The set $\cH_m^{\bar{\cC}}$, defined in \eqref{HmC}, will project to the hypersurface
\beq
\cF_m^{\bar{\cC}} := \{q\in M_0 \mid \|q_{C_1}\|\ge1, q_{C_2}\in S_m(q_{C_1}) \} 
\Leq{F:m:bar:C}
of the $(3d)$--dimensional center of mass configuration space 
\[M_0:=\{q\in M\mid q_N=0\}.\]
We choose one of the two continuous unit normal vector fields 
\beq
\textstyle 
N: \cF_m^{\bar{\cC}} \to T_{\cF_m^{\bar{\cC}}} M_0\qmbox{,} 
N(q)\in \operatorname{span}(\nabla \big\langle r-q_{C_1},\frac{q_{C_1}}{\|q_{C_1}\|}\big\rangle).
\Leq{def:N}
The normal component of $x\in \bR^d$, w.r.t.\  $\operatorname{span}(q_{C_1})$ equals 
\[ \textstyle  x^\perp:= x - \frac{\langle x,q_{C_1}\rangle}{\|q_{C_1}\|^2} q_{C_1}.\] 
In \eqref{the:surfaces} we now set (for a given parameter $ {\frak L}>0$) 
\begin{align}
\textstyle \cH_m^{\bar{\cC}}  := \{(p,q)  \in &\, \Sigma_{E,0} \mid
\, \|p_{C_2}\|\ge m , 
\|q_{C_1}\| \ge 1, q_{C_2}\in S_m(q_{C_1}), p(N(q)) < 0,  \NN \\
&  \|q^\perp_{C_2}\| \|p_{C_2}\| \le  {\frak L}, 
   \|p^\perp_{C_2}\| \|q_{C_2}\| \le {\frak L}, \|p^\perp_{C_1}\| \|q_{C_2}\| \le  {\frak L}\}.
\label{HmC}   
\end{align}
Physically, the condition $p(N(q)) < 0$ (with the pairing between $N(q)\in T_qM_0$ and  $p\in T^*_qM_0$)
means that the messenger cluster $C_2$ moves away from $C_1$, in the direction of $C_3$.
Mathematically, it implies that $I_m $ is a volume form on $\cH_m^{\bar{\cC}}$.
\begin{lemma}[Assumptions A1 and A2] \label{lem:A1:A2}
For all ${\frak L}>0$ the integrals
\[\textstyle 
I_m := \sum_{\bar{\cC}\in {\cal Q}} I_m^{\bar{\cC}}
\qmbox{with \footnotemark}  
I_m^{\bar{\cC}} := \int_{\cH_m^{\bar{\cC}}} {\cal V}_m\]
\footnotetext{for the pull-back ${\cal V}_m= {\iota}_m^* {\cal V}$ 
of ${\cal V}=\boldsymbol{i}_{X_H}\sigma_E$
w.r.t.\ the imbedding ${\iota}_m:\cH_m^{\bar{\cC}}\to \Sigma_{E,0}$ }%
are finite, $\,\lim_{m\to\infty}I_m=0$, and
the $(6d-2)$--form ${\cal V}$ on $\Sigma_{E,0}$ equals $\frac{\omega_0^{\,\wedge 3d-1}}{(3d-1)!}$.
\end{lemma}
The proof of that lemma will use a projection of the Poincar\'e surface $\cH_m^{\bar{\cC}}$ 
to the cotangent bundle of $\cF_m^{\bar{\cC}}$:
\begin{lemma} [Evaluation of the integrals] \label{lem:ev:int}
For $m\in \bN$ large 
\beq
I_m^{\bar{\cC}} \le
\int_{\cF_m^{\bar{\cC}}}\Big(
\int_{B_{r(q)}^{3d-1}}  \idty_{y(q)}(p) \,  dp \Big)d\cF(q),
\Leq{I:m:basic}
with $p \equiv (p_{C_1}, p^\perp_{C_2}, p^I_D)$, $q\equiv (q_{C_1}, q^\perp_{C_2}, q^I_{D})$, 
the Riemannian volume element $d\cF$ on $\cF_m^{\bar{\cC}}$ (w.r.t.\ the norm $\|\cdot\|_{\cM}$ on the 
tangent space, see \eqref{inner:product}), 
the ball $B_{r(q)}^{3d-1}$ of radius $r(q):=2(E-V(q))_+$ (w.r.t.\ $\|\cdot\|_{\cM^{-1}}$ on the 
cotangent space) and 
\begin{align*}
y(q) := \{p\in B_{r(q)}^{3d-1} \mid\, & \|p_{C_2}\|\ge m , 
\\
&  \|q^\perp_{C_2}\| \|p_{C_2}\| \le  {\frak L},
\|p^\perp_{C_2}\| \|q_{C_2}\| \le {\frak L}, \|p^\perp_{C_1}\| \|q_{C_2}\| \le  {\frak L}\}.\NN
\end{align*}
\end{lemma}
\begin{remark}[Independent variables]\quad\\ 
Note that in the definition of $y(q)$, 
$\|q_{C_2}\|$ and $\|p_{C_2}\|$ (unlike $\|q^\perp_{C_2}\|$ and $\|p^\perp_{C_2}\|$)
appear as functions of the integration variables, given by \eqref{F:m:bar:C}, respectively by
\[\textstyle
\frac{\|p_{C_2}\|^2}{2m_{C_2}}=E-V(q)-\frac{\|p_{C_1}\|^2}{2m_{C_1}}-\frac{\|p^I_D\|^2}{2m^I_{D}}.
\hfill \tag*{$\Diamond$}\]
\end{remark}
\textbf{Proof of Lemma \ref{lem:ev:int}:}\\
$\cH_m^{\bar{\cC}}$ projects diffeomorphically to its image in $T^*\cF_m^{\bar{\cC}}\subseteq T^*M_0$, via
\[n: \cH_m^{\bar{\cC}} \to T^*\cF_m^{\bar{\cC}} 
\quad\mbox{,}\quad(p,q)\mapsto \big(p-p(N(q),q) N^{\flat}(q)\big).\]
The cotangent bundle $T^*\cF_m^{\bar{\cC}}$ carries the canonical symplectic form $\omega_\cF$.
By Theorem C of \cite{FK1}, for the embedding $\imath:\cH_m^{\bar{\cC}} \to T^*M_0$
one has 
$\imath^*\omega_0 = n^*\omega_F$.
So 
\[\textstyle
\int_{\cH_m^{\bar{\cC}}} \frac{(\imath^*\omega_0)^{\wedge 3d-1}}{(3d-1)!} = 
\int_{n(\cH_m^{\bar{\cC}})} \frac{\omega_F^{\,\wedge 3d-1}}{(3d-1)!} \, ,\]
and we are left to determine the image $n(\cH_m^{\bar{\cC}})$, that is, $p-p(N(q))$, see \eqref{def:N}.

At the points where $q_{C_2}$ is linear dependent of  $q_{C_1}$, $q_{C_2}^\perp=0$ in
\beq
\textstyle
\nabla \big\langle q_{C_2}-q_{C_1},\frac{q_{C_1}}{\|q_{C_1}\|}\big\rangle 
=\sum_{i=1}^d 
\big(m^{-1}\frac{q_{C_1,i}}{\|q_{C_1}\|} + \langle q_{C_2}^\perp, f_i(q_{C_1})\rangle\big) \frac{\pa}{\pa q_{C_1,i}} 
+ \frac{q_{C_1,i}}{\|q_{C_1}\|}\frac{\pa}{\pa q_{C_2,i}} ,  
\Leq{projection:direction}
and thus its norm equals $\sqrt{1+1/m^2}$.
By comparison with the kinematical model, there is a constant such that after the $k$-th near-collision
$\|p_{C_1}\| \ge c_{\bar{\mu}} \bar{\mu}^k$, for any $\bar{\mu}\in(1,\mu^{1/2})$, 
with $\mu=  1+\big(\frac{m_{\min}}{m_{\max}} \big)^2$, see \eqref{en:expo}.
So by $\|q^\perp_{C_2}\| \|p_{C_2}\| \le  {\frak L}$, $\|q^\perp_{C_2}\|$ is exponentially small in $k$.
Dropping the normalizing factor $1/\sqrt{1+1/m^2}$ $\in(0,1)$, we can use \eqref{projection:direction} 
to estimate the integral.
\hfill $\Box$\\[2mm]
\textbf{Proof of Lemma \ref{lem:A1:A2}}\\
For an energy surface ${\iota}_E:\Sigma_E=H^{-1}(E)\to  T^*M^k$, 
equality of ${\cal V}=\boldsymbol{i}_{X_H}\sigma_E$ for the Liouville volume form $\sigma_E$ 
on $\Sigma_E$ with  $\frac{(\imath_E^*\omega_0)^{\wedge k-1}}{(k-1)!}$ was proven in \cite[(6.1)]{FK1}.

\noindent
Although \eqref{HmC} contains no explicit conditions concerning the nontrivial cluster $D\in \{C_1,C_2,C_3\}$, 
on $\cH_m^{\bar{\cC}}$ one has $\|p_{C_2}\|\ge m$ and thus
$H_D^I\le E- 2Cm^2+\cO(m^{-\alpha})$.
Henceforth we only consider\,\footnote{in particular $m\ge2$, so that $\|q_{C_2}\|$ is bounded away from zero,
see \eqref{def:hyper}.}
$m > m_0$, with $m_0\in \bN$ large enough such that 
$H_D^I\le -Cm^2$ (and that \eqref{I:m:basic} applies). This entails $\|q^I_D\|\le C m^{-2/\alpha}$.

After performing the spherical integrations for $p^\perp_{C_1}$, $p^\perp_{C_2}$ and $p^I_D$, 
given the norms 
\[ P \equiv (P_{C_1}, P^\perp_{C_1}, P^\perp_{C_2}, P^I_D )  := 
( \|p_{C_1}\|, \|p^\perp_{C_1}\|, \|p^\perp_{C_2}\|, \|p^I_D\| ) \] 
of the momenta, the r.h.s.\ of \eqref{I:m:basic} equals
\beq
c_1\int_{\cF_m^{\bar{\cC}}}\!\!\Big(\!
\int_{H_{r(q)}^{(4)}\cap Y(Q)}\!\! (P^\perp_{C_1})^{d-2} (P^\perp_{C_2})^{d-2} (P^I_D)^{d-1}\, 
d^4\!P \Big) d\cF(q)
\Leq{I:m:radial}
with 
\[\textstyle
H_{r(q)}^{(4)} := \Big\{P\in (\bR^+)^4\ \Big| \
\frac{P_{C_1}^2}{2m_{C_1}} + \frac{(P^I_{D})^2}{2m_D} \le r^2(q) - \frac{m^2}{2m_{C_1}}, 
P^\perp_{C_1}\le P_{C_1} , P^\perp_{C_2}\le P_{C_2} \Big\} , \]
the Riemannian measure $s_k$ of the sphere $S^k$, 
$c_1 := s_0^2 s_{d-2}^2 s_{d-1} $ and, for 
\[Q\equiv (Q_{C_1},Q^\perp_{C_1}, Q^\perp_{C_2} , Q^I_D) := 
(\|q_{C_1}\|, \|q^\perp_{C_1}\| , \|q^\perp_{C_2}\| , \|q^I_D\|),\]
\[Y(Q):=\{P\in (\bR^+)^4 \mid  
Q^\perp_{C_2} \le  {\frak L}/P_{C_2} ,
P^\perp_{C_2} \le {\frak L}/Q_{C_2}, P^\perp_{C_1} \le  {\frak L}/Q_{C_2} \}.\]
To perform the spherical integrations for the position $q$, we use that
for $m$ large, on $\cH_m$ the potential is dominated by $V(q)\le -2I (Q^I_{D})^{-\alpha}$, since
\begin{enumerate}[$\bullet$]
\item 
$\|q^I_D\|=\cO(m^{-2/\alpha})\ll m^{-1} =\operatorname{dist}\big(S_m(q_{C_1}),q_{C_1}\big)$, and
\item 
by the condition $q_N=0$, we have 
$q_{C_3} = 
- \frac{m_{C_1}+m_{C_2}(1-1/(m \|q_{C_1}\|))}{m_{C_3}}q_{C_1} 
- \frac{m_{C_2}}{m_{C_3}}q_{C_2^\perp}$ 
so that $\|q_{C_3}-q_{C_i}\| \ge \frac{m_{\min}}{m_{\max}}$ for $i=1,2$, using $\|q_{C_1}\|\ge1$.
\end{enumerate}
So, with $R(Q):= 2(E+2I (Q^I_D)^{-\alpha})$ and 
$F_m := \{Q\in (\bR+)^4 \mid Q_{C_1}\ge 1\}$, \eqref{I:m:radial} is dominated for $d>2$ by
\begin{align}
&c_2 \int_{F_m}\!\!\Big(\!
\int_{H_{R(Q)}^{(4)}\cap Y(Q)}\!\! (P^\perp_{C_1})^{d-2} (P^\perp_{C_2})^{d-2} (P^I_D)^{d-1}\, 
d^4\!P \Big) Q_{C_1}^{d-1} (Q^\perp_{C_2})^{d-2} (Q^I_D)^{d-1}\, d^4Q \NN\\
&\le 
c_3\!\! \int_1^\infty\!\!\!\!  \int_0^{m^{-\alpha}}\!\!\!\!\!\!\!\Big(\!
\int_m^\infty\!\!\! \int_0^{c (Q^I_{D})^{-\alpha}}\hspace{-10mm}(P_{C_2})^{1-d}(P^I_D)^{d-1} 
dP_{C_2} \, dP^I_D \Big) Q_{C_2}^{2-2d} Q_{C_1}^{d-1}  (Q^I_D)^{d-1}\, dQ^I_D\,dQ_{C_1} \NN\\
&\le 
c_4 \, m^{2-d}\,\int_1^\infty Q_{C_1}^{1-d} \,dQ_{C_1}  \int_0^{m^{-\alpha}}\!
  (Q^I_D)^{d(1-\alpha/2)-1}\, dQ^I_D 
\ \le \
c_5 \, m^{2-d-\alpha d(1-\alpha/2)}
\label{I:m:radial:final}
\end{align}
Thus for $d\ge 3$ dimensions the integrals $I_m$ are finite, and $\lim_{m\to\infty}I_m=0$.
\hfill$\Box$\\[-3mm]
\begin{remark}[The case of two dimensions]\quad\\ 
Note that although the integral $\int_1^\infty Q_{C_1}^{1-d} \,dQ_{C_1}$ diverges logarithmically for
$d=2$, the power of $m$ in \eqref{I:m:radial:final} is negative also in this case.
 
Probably one could extend our main result (page \pageref{thm:main}) to $d\ge2$
by taking into account that the messenger particle has to
return to $C'_1$ after having experienced its near-collision with the cluster $C_3$.
However, to prove that this effect leads to an additional negative power of  $Q_{C_1}$
in \eqref{I:m:radial:final} would require additional work.
\hfill $\Diamond$
\end{remark}
\begin{lemma}  
If $x\in\operatorname{NAV}$, then for ${\frak L}$ large enough in the definition \eqref{the:surfaces}
of Poincar\'e surfaces, the forward orbit hits almost all $\cH_m$.  
\end{lemma}
\textbf{Proof:}\\
We check the conditions in Def.\ \eqref{HmC} of the Poincar\'e surfaces $\cH_m^{\bar{\cC}}$,
$\bar{\cC}\in {\cal Q}$ :
\begin{enumerate}[$\bullet$]
\item
$\|\tilde{p}_{C_2}(t)\|\ge m$ holds for $t<T(x)$ large by  Lemma \ref{lem:part:energy}, in combination 
with Theorem \ref{thm:NAV:wandering}.
\item 
$\|\tilde{q}_{C_1}\|(t) \ge 1$ follows for $t<T(x)$ large by von Zeipel's Theorem \ref{thm:von:Zeipel}.
\item
$\tilde{q}_{C_2}(t)\in S_m(\tilde{q}_{C_1}(t))$ with $\tilde{p}(N(\tilde{q})) < 0$ holds for a sequence
of times converging to $T(x)$, since, as shown in Section \ref{sub:sect:B2}, the messenger cluster $C_2$
moves infinitely many times between disjoint neighborhoods of  $\tilde{q}_{C_1}$ and $\tilde{q}_{C_3}$.
As proven in Lemma \ref{lem:three:particles:int} (see Assertion 2. of Lemma \ref{lem:three:particles:no:int} for
the (converse) estimate)), these neighborhoods are of radius 
$\cO\big(\big({\tilde{K}_{C_1,C_2} }\big)^{-1/\alpha}\big) =
\cO\big(\big({\tilde{K}_{C_2} }\big)^{-1/\alpha}\big)$ which shrinks to zero as $t\nearrow T(x)$
by Lemma \ref{lem:part:energy}. So for any given $m\in \bN$ they lie in different half-planes
of $\bR^d$, defined by $S_m$, see \eqref{def:hyper}.
\item 
The following three statements \eqref{three:statements}
all hold for parameter ${\frak L}\ge 2\|\tilde{L}(0)\|$ and $t<T(x)$ large,
since
\begin{enumerate}[(a)]
\item 
By Lemma \ref{lem:small:var:a:m}.3 the deviations of the trajectories $t\mapsto \tilde{q}_{C_i}$ 
of the cluster centers from motion on affine lines 
of $\bR^d$ (whose definition is given in the lemma) becomes small.
\item 
By Lemma \ref{lem:three:particles:int} $C_2$ has close encounters with $C_1$ and $C_3$, see 
above.
\item 
At these Poincar\'e times $t$ the nontrivial cluster $D(t)$ has a size that goes to zero 
(see Assertion 2. of Lemma \ref{lem:three:particles:no:int}) and thus
is compared to the minimal distance $\cO(1/m)$ of the cluster centers.
\end{enumerate}
This is reflected in the angular momenta:\\ 
By Lemma \ref{lem:cluster:ang:mom} the cluster angular
momenta $\tilde{L}^E_C(t) \| \le \|L(x)\| + \ell$ for all $\ell>0$ and $t\in[t_0, T(x))$, $C\in\cA(t)$ 
for $t_0(\ell)\in (0,T(x))$. Similarly, $\|\tilde{L}^I_C(t)\| \le \ell$.

By Lemma \ref{lem:small:var:a:m}.1, the variations of the relative angular momenta go to zero.

Thus we can compare with the straight line geometry of the kinematical model to see that
for the above value of ${\frak L}$ the violation of one of the three conditions
\beq
\|q^\perp_{C_2}\| \|p_{C_2}\| \le  {\frak L}, 
   \|p^\perp_{C_2}\| \|q_{C_2}\| \le {\frak L}, \|p^\perp_{C_1}\| \|q_{C_2}\| \le  {\frak L}
 \Leq{three:statements}
would be in contradiction with $x\in \operatorname{NAV}$ and $t<T(x)$ large.
\hfill $\Box$
\end{enumerate}

%
\addcontentsline{toc}{section}{References}
%

\end{document}